\documentclass[11pt]{article}

\usepackage{color}
\usepackage{amsmath}
\usepackage{mathrsfs}
\usepackage{amssymb}
\usepackage{amsthm}
\usepackage{epstopdf}
\usepackage{graphicx}
\usepackage{mathtools}
\usepackage[usenames,dvipsnames, table]{xcolor}
\usepackage[hang]{caption}
\usepackage{hyperref}
\usepackage{cite}

\usepackage{oldgerm}  
\usepackage[yyyymmdd,hhmmss]{datetime}

\usepackage{mathabx} 
\def\bee{\begin{enumerate}}\def\eee{\end{enumerate}}
\def\bei{\begin{itemize}}\def\eei{\end{itemize}}
\newcommand{\nco}{\newcommand}

\nco{\red}{\color{red}}
\nco{\blue}{\color{blue}}
\nco{\cyan}{\color{cyan}}
\nco{\brown}{\color{Magenta}}

\nco{\magenta}{\color{magenta}}

\nco{\violet}{\color{violet}}
\nco{\olive}{\color{Emerald}}
\nco{\orange}{\color{orange}}
\nco{\redend}{\normalcolor}
\nco{\blueend}{\normalcolor}
\def\inv#1{\frac{1}{#1}}
\def\tr{{\rm tr}\,}

\def\dot{point}
\def\eg{{\it e.g.\,}}
\def\ommit#1{{}}
\def\({\left(}
\def\){\right)}
\def\ie{{\it i.e.,\/}\ }
\def\ie{{\rm i.e.,\/}\ }

\definecolor{cb}{rgb}{.8,.5,0}

\nco{\rnc}{\renewcommand}
\rnc{\title}[1]{{\Large\bf\mbox{}\\\medskip#1\bigskip\medskip\\}}
\rnc{\author}[1]{{\large #1\smallskip\\}}
\nco{\address}[1]{{\em #1\medskip\\}}

\def\be{\begin{equation}}\def\ee{\end{equation}}
\def\bea{\begin{eqnarray}}\def\eea{\end{eqnarray}}
\def\bee{\begin{enumerate}}\def\eee{\end{enumerate}}
\def\bei{\begin{itemize}}\def\eei{\end{itemize}}
\def\oh{\frac{1}{2}}
\def\ommit#1{{}}
\def\omit#1{{}}
\def\CO{{\mathcal O}}\def\CP{{\mathcal P}}\def\CS{{\mathcal S}}
\def\Cs{{\mathfrak s}}
\def\k{\kappa}
\def\Ge{\epsilon}
\def\tx{\tilde x}
\def\tX{\tilde X}
\def\X1{\X1}
\def\X1{Y_2}
\def\tX1{\tilde \X1}
\def\tX2{\tilde X_2}
\def\tX{\tilde X}
\def\tY{\tilde Y}
\def\Arcsin{\mathrm{Arcsin\,}}
\def\Arccos{\mathrm{Arccos\,}}
\def\[{[\![} \def\]{]\!]}

\def\comment#1{{\blue #1}}
\def\comment#1{}
\def\adots{\mathinner{\mkern2mu\raise1pt\hbox{.}\mkern3mu\raise4pt\hbox{.}\mkern1mu\raise7pt\hbox{.}}}  
\def\dddots{\mathinner{\mkern2mu\raise9pt\hbox{.}\mkern3mu\raise3pt\hbox{.}\mkern2mu\raise-4pt\hbox{.}}}

\def\new#1{{\magenta #1}}
\def\new#1{#1}
\def\neww#1{{\cyan #1}}
\def\neww#1{#1}
\def\newa#1{{\cyan #1}}
\def\newa#1{#1}
\def\moi#1{{#1}}\def\moi#1{}

\begin{document}
\begin{titlepage}
%
\begin{center}
\title{Counting partitions by genus\\[8pt]
 I. Genus 0 to 2}
\medskip
\author{Jean-Bernard Zuber}
\address{Sorbonne Universit\'e,  UMR 7589, LPTHE, F-75005,  Paris, France\\ \& CNRS, UMR 7589, LPTHE, F-75005, Paris, France\\
{\tt zuber@lpthe.jussieu.fr}}
\bigskip\bigskip
\begin{abstract}
The counting of partitions according to their genus 
is revisited. The case of genus 0 --non-crossing partitions-- is well known.  
Our approach relies on two pillars : first a functional equation between 
generating functions, originally written in genus 0 and interpreted graphically by Cvitanovic, is generalized to higher genus; 
secondly, we  show that all partitions may be reconstructed  from the ``(semi)-primitive" ones
introduced  by Cori and Hetyei. 
  Explicit results  for the generating functions of all types of partitions are obtained in genus 1 and 2.
  \new{This gives a second order interpolation between expansions on ordinary or on free cumulants.}
\end{abstract}
\end{center}

\bigskip

\normalcolor
 \end{titlepage}

\bigskip\bigskip
\section{Introduction}

\subsection{Partitions, their genus and their census}
\omit{\cyan We are interested in partitions of the {\it set} $\[n\]:=\{1,\cdots,n\}$. 
If $\alpha$ is such a partition, made of $\alpha_1$ parts of length 1,  $\alpha_2$ parts of length 2, etc, 
we write  $\alpha\in \CP(n)$ and say that $\alpha$ is of {\it type}  $[\alpha]=[1^{\alpha_1}, \cdots, n^{\alpha_n}]$,
which may be regarded as a partition of the {\it integer} $n$: $[\alpha]\vdash n$.\\}
{Consider the set $\CP(n)$ of partitions of the {\it set} $\[n\]:=\{1,\cdots,n\}$. 
If $\alpha\in \CP(n)$ is made of $\alpha_1$ parts of length 1,  $\alpha_2$ parts of length 2, etc, 
we  say that $\alpha$ is of {\it type}  $[\alpha]=[1^{\alpha_1}, \cdots, n^{\alpha_n}]$,
which may be regarded as a partition of the {\it integer} $n$: $[\alpha]\vdash n$.\\}
\omit{\cyan Note that when listing the parts of a partition $\alpha= (\{i_1\}, \cdots \{i_{\alpha_1}\},\{j_1,j_2\}, \cdots) $,\\
(i) the ordering of elements in each part is immaterial, and we thus choose to write them in increasing order;\\
(ii) the relative position of parts is immaterial. }

Let  $C_{n,[\alpha]}$ denote 
 the number of partitions of type $[\alpha]$ 
 \be\label{sumovergenus}C_{n,[\alpha]}=\frac{n!}{\prod_{\ell=1}^n \alpha_\ell!  (\ell!)^{\alpha_\ell}}\,. \ee

The census of partitions may be subject to different conditions. 
 In particular, it is well known, as we recall below in sect.~\ref{genus-def}, 
 that any partition  $\alpha$ may be assigned a {\it genus} $g(\alpha)$ by a formula descending from 
Euler's relation. Curiously, the census of partitions according to their genus is still an open problem, in spite of several
fundamental contributions, \cite{Krew, WL1, WL2, CoriH13, CoriH17}. Except for a few particular cases, only 
the case of genus 0 is thoroughly known: the {\it non crossing partitions} (or planar) have  been enumerated 
by Kreweras \cite{Krew}, before reappearing in various contexts, matrix integrals \cite{BIPZ, Cvitanovic}, 
free probability \cite{V86, Speicher}, and more recently, out-of-equilibrium quantum systems~\cite{HB22,FKP22}.
\\
\omit{\olive Now, it is also well known, as we recall below in sect.~\ref{genus-def},  that any partition  $\alpha\in \CP(n)$ may be assigned a {\it genus} $g(\alpha)$ by a formula descending from 
Euler's relation. It is then natural  to refine the counting of partitions according to their genus.}\\
Let $C^{(g)}_{n,[\alpha]}$ denote the number of partitions in $\CP(n)$ of type $[\alpha]$ and genus $g$. Obviously $\sum_g C^{(g)}_{n,[\alpha]}=C_{n,[\alpha]}$.

We find it convenient to use generating functions (GF) to encode these numbers.
Introduce a set of indeterminates $\kappa_n$, $n\in \Bbb{N}_+$, their GF
\be\label{defW} W(x)=\sum_{n\ge 1} \kappa_n x^n \ee
and then
 \bea
 \label{defZg}  Z(x)&=& 1+\sum_{n\ge 1}\sum_{[\alpha]\vdash n} C_{n,[\alpha]} \kappa_{[\alpha]} x^n
 =\sum_g    Z^{(g)}(x)\, \\
 \nonumber 
  Z^{(g)}(x)&=& \delta_{g0}+\sum_{n\ge 1}\sum_{[\alpha]\vdash n} C^{(g)}_{n,[\alpha]} \kappa_{[\alpha]} x^n\, \eea
where 
\be\label{falpha}\k_{[\alpha]} := \prod_{\ell=1}^n   \k_\ell^{\alpha_\ell}
\,. \ee

There is a well known relation between $Z^{(0)}$  and $W$, which has been found in different
avatars \cite{BIPZ, Cvitanovic,Speicher}, see below (\ref{Z0W}).\\
 To extend such a relation to higher genus, we rely on a proven method. The diagrams encoding 
the partitions are first reduced to basic diagrams, in finite number at given genus. In a second step, 
all diagrams --all partitions-- are reconstructed by ``dressing" the basic ones. This method is 
well known in combinatorics and quantum field theory (``skeleton diagrams"). 
In the context of enumeration of unicellular maps and partitions of given genus, 
it has been explored by Chapuy~\cite{Chapuy} and Cori--Hetyei~\cite{CoriH17}, who call 
the basic diagrams {\it schemes} and {\it primitive}, respectively.\\
In this paper, explicit formulae relating $Z^{(g)}$  to $W$ and its derivatives will be found  for $g=1,2$, 
and the corresponding expressions of $Z^{(g)}(x)$ are given by Theorem 1, 
     (\ref{Zgenus1}),
     and Theorem 2,  (\ref{Z2genus2}).  Extension to higher genera is in principle feasible, if the list of their primitive diagrams        is known.\\[6pt]
 
 \subsection{Genus dependent cumulant expansion}    
     The question of the relation between $Z$ and $W$ 
     also arises in probability theory and statistical mechanics. There,  it is common practice to associate 
{\it cumulants} to moments of random variables. 
If $X$ is a random variable with moments $m_n=\Bbb{E}(X^n)$ of arbitrary order, we decompose these moments on 
cumulants $\k_m$ and their products associated with {partitions}   $\alpha\in \CP(n)$  
\be\label{defcumul0} m_n=\sum_{\alpha\in \CP(n)}\k_\alpha\,.\ee
Thus each term in (\ref{defcumul0}) may be regarded as associated with a splitting  of $\[n\]$ into parts described  by the partition:
in statistical mechanics, the terms $\k_\alpha$ are dubbed the {\it connected parts} of the moment $m_n$. 
If we assume  that 
the moments and cumulants do not depend on extra variables like momenta, etc, we may rewrite (\ref{defcumul0})
as
\be\label{defcumul} m_n=\sum_{[\alpha]\vdash n} C_{n,[\alpha]}\k_{[\alpha]}\,.\ee
For example, $m_4=   \k_4+  4\,\k_3  \,\k_1+ 3 \k_2^2   + 6\,\k_2  \,{\k_1}^2 + {\k_1}^4\,.$ \\
Thus the indeterminates $\k_\ell$  have acquired the  meaning of cumulants, and
$Z(x)$ and $W(x)$ are the GF of moments and cumulants, respectively. \\

Then,  making use of the genus $g(\alpha)$ mentioned above, 
it is  natural to modify the expansion  (\ref{defcumul}) by weighting the various terms acccording to their genus.
Introducing a parameter $\epsilon$, we write
\be \label{cumuleps}  m_n(\epsilon) = \sum_{\alpha\in \CP(n)} \epsilon^{g(\alpha)}\k_\alpha\ee 
or 
\be \label{cumuleps}  m_n(\epsilon) = \sum_{[\alpha]\vdash n} \sum_{g=0}^{g_{{\mathrm max}}([\alpha])} C^{(g)}_{n,[\alpha]}\epsilon^{g}\k_{[\alpha]}\\,.\ee 
\\
For example, $m_4(\epsilon)=   \k_4+  4\,\k_3  \,\k_1+ (2+\epsilon) \k_2^2   + 6\,\k_2  \,{\k_1}^2 + {\k_1}^4\,,$ {see below}.

Obviously for $\epsilon=1$, we 
recover the usual expansion  (\ref{defcumul}), whereas for $\epsilon=0$, we have an expansion on {\it non crossing} (or free, or planar)
cumulants.
Thus (\ref{cumuleps}) provides an interpolation between the usual cumulant expansion and that on non crossing ones.

\omit{In this paper, we try to determine the numbers $ C^{(g)}_{n,[\alpha]}$.  Or alternatively, we strive to find relations between the
(ordinary) generating functions  (GF) of the $m_n(\Ge)$ and of the \neww{$\k_\ell$}:
\bea \label{defZj}  Z(x, \epsilon)&=& 1+ \sum_{n\ge1} m_n(\epsilon) x^n \\
    \nonumber &=& \sum_{g\ge 0}  Z^{(g)}(x) \epsilon^g \\
    \label{defWj} W(x)&=& \neww{\sum_{\ell\ge 1} \kappa_\ell x^\ell }\,.
    \eea
    This will be achieved for genus 1 and 2, and the corresponding expressions of $Z^{(g)}(x)$ are given by Theorem 1, 
     (\ref{Zgenus1}),
     and Theorem 2,  (\ref{Z2genus2}).  Extension to higher genera is in principle feasible, if the list of their primitive diagrams 
     is known.}
    

\subsection{Eliminating or reinserting singletons}
In a partition, parts of size 1 are called {\it singletons}. It is natural and easy to remove them in the counting, or
 to relate the countings of partitions with or without singletons. Let us denote with a hat the GF of partitions without 
singletons: $\hat  Z^{(g)}(x)$,  and derive the relation between $\hat  Z^{(g)}(x)$ and $Z^{(g)}(x)$.
This is particularly easy in the language of statistics, where discarding singletons amounts to going to
a centered variable: $X=\hat X +\Bbb{E}(X)= \hat X+m_1=\hat X+\k_1$
$$ m_n=\Bbb{E}(X^n) = \Bbb{E}((\hat X+\k_1)^n)= \sum_{r=0}^n {n\choose r} \, \hat m_{n-r} \,\k_1^k$$
and, since singletons do not affect the genus, see below
sect.~\ref{reduction},
\be\label{wwosingl} C^{(g)}_{n,[\alpha', 1^r]} ={n\choose r} C^{(g)}_{n-r, [\alpha']}\, \ee
where the partition $\alpha'$ is {\it singleton free} (s.f.). 
For example, 
\bea \nonumber  \hskip-20mm m_1\!\!&=&\!\! \k_1 \\ \nonumber
m_2 \!\!&=&\!\! \k_2+\k_1^2 \\ \nonumber
m_3 \!\!&=&\!\! \k_3+ 3 \k_2\k_1+\k_1^3 \\ \nonumber
m_4 \!\!&=&\!\! \k_4+ (2+\epsilon) \k_2^2 +4 \k_3\k_1 +6\k_2\k_1^2  +\k_1^4\\
\nonumber  m_5\!\!&=&\!\! {\k_5 
 +  5\,\k_4\,\k_1+5 (1+\Ge) \k_3 \k_2  
+ 10\,\k_3 \,{\k_1}^2+ 5 (2+\Ge)  \,{\k_2}^2 \k_1  + 10\,\k_2  \,{\k_1}^3 
 + {\k_1}^5  \,, }\eea
 etc.

 Then
\bea \nonumber  Z^{(g)}(x) &=& \sum_{n\ge 0}  x^n\sum_{[\alpha]\atop \alpha\in \CP(n)}  C^{(g)}_{n,[\alpha]} \k_{[\alpha]}\\
\nonumber  &=& \sum_{n\ge 0}  x^n  \sum_{r=0}^n  \sum_{[\alpha']\atop \alpha' \in \CP(n-r),\, \mathrm{s.f.}} C^{(g)}_{n,[1^r,\alpha']}   \k_{[\alpha']}  \k_1^r \\
\nonumber  &=&  \sum_{n'\ge 0}  x^{n'}\sum_{[\alpha']\atop  \alpha' \in \CP(n'),\,\mathrm{s.f.}} C^{(g)}_{n',[\alpha']}   \k_{[\alpha']}  \sum_{r\ge 0}  {n'+r\choose r}   \k_1^r x^r\\
\nonumber  &=&  \sum_{n'\ge 0}  \sum_{[\alpha']\atop  \alpha' \in \CP(n'),\,\mathrm{s.f.}} C^{(g)}_{n',[\alpha']}   \k_{[\alpha']}  \frac{x^{n'}}{(1-\k_1 x)^{n'+1}}\\
\label{ZZhat} &=&  \inv{1-\k_1 x}  \, \hat  Z^{(g)}\big(\frac{x}{1-\k_1 x}\big) \,,
\eea
 and conversely
\be\label{ZhatZ} \hat Z^{(g)} (u) = \inv{1+\kappa_1 u}\, Z^{(g)}\big(\frac{u}{1+  \kappa_1 u}\big)\,.\ee

 \omit{Reinstating the singletons is then an easy matter. 
The expression of $C^{(g)}_{n,\alpha}$ for $\alpha=[\alpha', 1^r]$ containing $r$ $1$'s is readily obtained from the graphical representation of
Fig. \ref{Fig1}. The insertion of the $r$ $1$'s can be done anywhere in the diagram of $\alpha'$, and does not affect the genus, thus
\be\label{insert1} C^{(g)}_{n,[\alpha', 1^r]} ={n\choose r} C^{(g)}_{n-r, [\alpha']}\,. \ee
Hence 
\bea \nonumber  \hskip-20mm m_2 \!\!&=&\!\! \k_2+\k_1^2 \\ \nonumber
m_3 \!\!&=&\!\! \k_3+ 3 \k_2\k_1+\k_1^3 \\ \nonumber
m_4 \!\!&=&\!\! \k_4+ (2+\epsilon) \k_2^2 +4 \k_3\k_1 +6\k_2\k_1^2  +\k_1^4\\
\nonumber  m_5\!\!&=&\!\! {\k_5 
 +  5\,\k_4\,\k_1+5 (1+\Ge) \k_3 \k_2  
+ 10\,\k_3 \,{\k_1}^2+ 5 (2+\Ge)  \,{\k_2}^2 \k_1  + 10\,\k_2  \,{\k_1}^3 
 + {\k_1}^5  \,, }\eea
 etc.
Making use of the relation  (\ref{insert1})  below, namely
$$ C^{(g)}_{n,[\alpha', 1^r]} ={n\choose r} C^{(g)}_{n-r, [\alpha']}\, $$
for $[\alpha']$ ``singleton free" (s.f. in short),
we may write for the genus $g$ GF  
\bea   \nonumber Z^{(g)}(x; \kappa_i)&=& \sum_n  x^n \sum_{[\alpha]\vdash n}  C^{(g)}_{n,[\alpha]} \kappa_{[\alpha]} \\
\nonumber&=& \sum_n x^n \sum_{[\alpha'] \,\mathrm{s.f.}} \sum_{r\le n} {n\choose r}  \kappa_1^r C^{(g)}_{n-r, [\alpha']} \kappa_{[\alpha']} \\
\nonumber&=& \sum_{n'\ge 0} \sum_{r\ge 0} {n'+r\choose r}  \kappa_1^r x^r   \sum_{[\alpha']\, \mathrm{s.f.}} C^{(g)}_{n', [\alpha']} \kappa_{[\alpha']} x^{n'}\\
\nonumber&=& \sum_{n'\ge 0} \frac{1}{(1-x \kappa_1)^{n'+1}}    \sum_{[\alpha'] \mathrm{s.f.}} C^{(g)}_{n', [\alpha']} \kappa_{[\alpha']} \\
\label{singlefree}&=& \frac{1}{(1-x \kappa_1)}  Z^{(g)\, \mathrm{s.f.}} \big(\frac{x}{1-x \kappa_1}\big)
\eea
where the last GF denotes the s.f. $Z^{(g)\, \mathrm{s.f.}} (u)\equiv  Z^{(g)}(u; \kappa_1=0, \kappa_{i\ge 2})$. Conversely
\be Z^{(g)\, \mathrm{s.f.}} (u) = (1+\kappa_1 u)\, Z^{(g)}\big(\frac{u}{1+ u \kappa_1}\big)\,.\ee
}


\section{Partitions and their genus}
In this section, we recall some standard notions on partitions, show how to associate a graphical 
representation to a partition and introduce its genus in a natural way. 

\subsection{Parts of a partition} 
\label{subsec21}
As explained in sect. 1, we are interested in partitions of the {\it set} $\[n\]$. 
\omit{\magenta If $\alpha\in \CP(n)$ is such a partition, made of $\alpha_1$ parts of length 1,  $\alpha_2$ parts of length 2, etc, 
we say that $\alpha$ is of type  $[\alpha]=[1^{\alpha_1}, \cdots, n^{\alpha_n}]$,
which may be regarded as a partition of the {\it integer} $n$: $[\alpha]\vdash n$.}\\
Note that  when listing the parts of a partition $\alpha= (\{i_1\},\cdots \{i_{\alpha_1}\},\{j_1,j_2\}. \cdots) $,\\ 
(i) the ordering of elements in each part is immaterial, and we thus choose to write them in increasing order;\\
(ii) the relative position of parts is immaterial.

For example, consider the partition $(\{1,3,4,6,7\},\{2,5,9\},\{8\},\{10\})$ of $\[10\]$. It is of type 
$[1^2,3,5]$ with two singletons $\{8\}$ and $\{10\}$.
Clearly the order of elements within each part is irrelevant, \eg  parts  $\{1,3,4,6,7\}$
and $\{3,4,1, 7,6\}$ describe the same subset of $\[10\]$. One may thus order the elements of each part. Likewise 
the relative order of the parts is immaterial: \\ $(\{1,3,4,6,7\},\{2,5,9\},\{8\},\{10\})$ and  $(\{2,5,9\},\{8\},\{1,3,4,6,7\},\{10\})$ describe the same partition. \\[4pt]

 \begin{figure}\begin{center}
{\includegraphics[width=.5\textwidth]{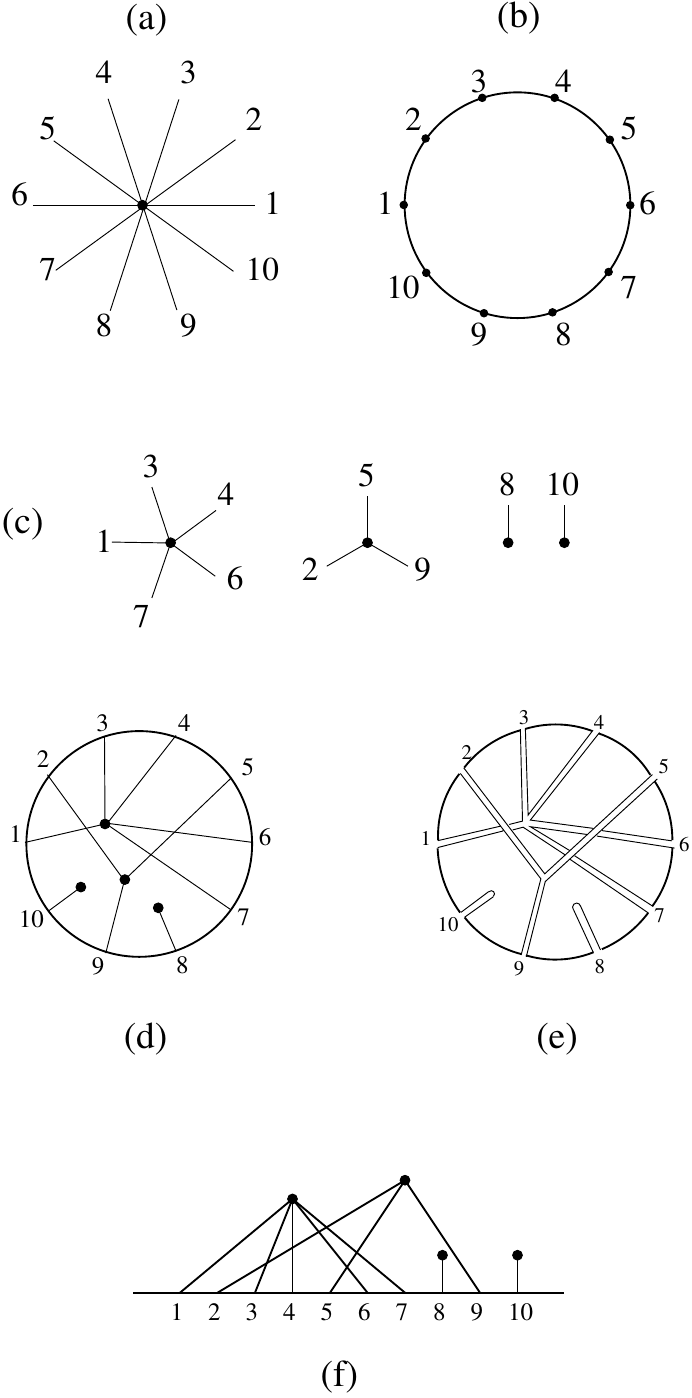}}
\end{center} 
\caption{\small{The  partition $(\{1,3,4,6,7\},\{2,5,9\},\{8\},\{10\})$ of $\[10\]$.  (a) and (b): two equivalent representations of the 10-vertex; (c) the four other vertices; (d) a   contribution to $C^{(g)}_{10, [1^2\,3 \, 5]}$; (e)  the double line version of (d), with  three faces and thus genus $g=$ 2;
(f) the linear version of (d).}}
\label{Fig1}
\end{figure}


\subsection{Combinatorial and graphical representations of a partition and its genus}
\label{genus-def}
 A general partition $\alpha$ of $\CP(n)$ may be described in terms of  
 a pair of {\it permutations} $\sigma$ and $\tau$, both in $\mathcal{S}_n$: $\sigma$ is the cyclic permutation
$(1,2,\cdots,n)$;  $\tau$ belongs to the class $[\alpha]$ of $\mathcal{S}_n$, and its {\it cycles} are 
described by the parts of $\alpha$, 
thus subject to the condition  (i) above:  {\it each cycle is  an increasing list of integers}.\\
The genus  $g$ of the partition  is then defined by \cite{Jacques}
\be n+2-2g=  \#\mathrm{cy}(\tau)+\#\mathrm{cy}(\sigma) + \#\mathrm{cy}(\sigma\circ \tau^{-1})\,\ee 
or in the present case, 
\be\label{genus} -2g = \sum \alpha_\ell -1 - n  + \#\mathrm{cy}(\sigma\circ \tau^{-1}) \,.\ee
since here $\#\mathrm{cy}(\sigma) =1$ and  $\#\mathrm{cy}(\tau)= \sum \alpha_k$.
Since $ \#\mathrm{cy}(\sigma\circ \tau^{-1})\ge1$, we find an upper bound on $g$
\be\label{gmax} g\le g_{\mathrm{max}}:= \bigg\lfloor\oh(n -\sum \alpha_k)\bigg\rfloor \,, \ee
see also \cite{Yip}.  We recall below why this definition of the genus is natural. 
\\[4pt]
Example.\\ For the above partition of $\[10\]$, $\sigma=(1,2,\cdots,10)$, $\tau=(1,3,4,6,7)(2,5,9)(8)(10)$,\\
$\sigma\circ\tau^{-1}=\new{(1,8,9,6,5,3,2,10)}(4)(7)$.  
Thus $2g=11-4-3=4$, $g=2$, while $g_{\mathrm{max}} =3$. \\[4pt]

To a given partition, we may also attach a map: it has  $\alpha_\ell$  $\ell$-valent vertices, \new{in short {\it $\ell$-vertices}}~\footnote{\label{footn}Remember that $\alpha_\ell$ are the multiplicities introduced
in (\ref{falpha})}, for $\ell=1,2,\cdots$,   whose edges are numbered clockwise by the elements of the partition, and
a special $n$-valent vertex, with its $n$ edges numbered {\it anti}-clockwise from 1 to $n$,
 see Fig.~\ref{Fig1}a,c. \newa{Edges are connected pairwise by matching their indices.
 Two maps are regarded as topologically equivalent if they encode the same partition.} In fact
it is topologically equivalent and more handy to attach  $n$ \dot s {\it clockwise}  on a circle, and to connect them pairwise by arcs of the circle,  
see Fig.~\ref{Fig1}b. 
Now the permutation $\sigma$ describes the connectivity of the $n$ \dot s on the circle, while $\tau$ describes how these \dot s are
connected through the other vertices. It is readily seen that the permutation $\sigma\circ\tau^{-1}$ describes the circuits bounding clockwise
the faces of the map. This is even more clearly seen if one adopts a double line notation for each edge \cite{tH}, thus transforming 
the map into a ``fat graph", see Fig.~\ref{Fig1}e . Thus the number of cycles of $\sigma\circ\tau^{-1}$  is the number $f$ of faces of the map. 
Since each face is homeomorphic to a disk,  gluing a disk to each face transforms the map into a closed Riemann surface, to
which we may apply  Euler's formula
\be\label{Euler}2 -2g = \#(\mathrm{vertices})-\#(\mathrm{edges})+\#(\mathrm{faces})= 1+\sum_\ell \alpha_\ell - n +f\ee
with $f=\#\mathrm{cy}(\sigma\circ \tau^{-1})$,  and we have reproduced (\ref{genus}). \\

Remark 1. This coding of a map, or here of a partition, by a pair of permutations, with a resulting expression of its genus,  is an old 
idea originating in the work of Jacques, Walsh and Lehman \cite{Jacques, WL1, WL2} 
and rediscovered and used with variants by many authors since then
\cite{Drouffe}. 

Remark 2. The diagrammatic representation that we adopt here differs from that of other authors \cite{Yip, CoriH17}: in fact it is the 
{\it dual} picture, with our vertices corresponding to faces of these authors. Our preference for the former is due to its analogy with Feynman diagrams\dots


\subsection{Glossary}
It may be useful to list some elements of terminology used below.\\
 It is convenient to represent a partition of $\CP(n)$ by a diagram. It may be a {\it circular}
 diagram, with $n$ points equidistributed clockwise, as on Fig. \ref{Fig1}-d, and it has a genus as explained above.
 \new{We distinguish the {\it points} on the circle from the {\it vertices} which lie inside the disk.}\\
 Occasionally we use a {\it linear} diagram, with $n$ points labelled from 1 to $n$ on a line (or an arc), \new{and vertices above the line.}\\ 
Note that if we give each point of the circle a weight $x$ and each $k$-vertex the weight $\kappa_k$, 
the sum of diagrams of genus $g$ builds the GF $Z^{(g)}(x)$.\\
In a (circular) diagram, we call {\it 2-line} a pair of edges attached to a 2-vertex. 
In the following, the middle 2-vertex will be omitted on 2-lines, 
to avoid overloading the figures. A 2-line is then just a straight line between two points of the circle. \\
In a diagram, we call {\it adjacent} a pair of edges joining a vertex to adjacent points on the circle.
\neww{For example, on Fig. \ref{Fig2}, the edges ending at 1 and 3 are not adjacent, those ending at 3 and 4 are.}\\
In the following discussion, it will be important to focus on a  point on the circle, say point 1, and see what it is connected to. We  shall
refer to it as the {\it marked point}.\\
 If $\alpha$ is a partition of $\CP(n)$ of a given type, all its conjugates by \new{powers of} the cyclic permutation $\sigma$ have
 the same type. Counting partitions of a given type thus amounts to counting {\it orbits} \new{of diagrams} 
 under the action of $\sigma$, while recording the length (cardinality) of each orbit.  
 Diagrammatically, the point 1 being \new{marked}, we list orbits under rotations of the inner pattern of vertices and edges by 
 the cyclic group $\Bbb{Z}_n$, and record the length of each orbit. An orbit has \new{a weight equal to its} length  $n/\Cs$, where $\Cs$ is the order of the \newa{stabilizer of the diagram -- a  subgroup of the rotation group}.
 For example, the left-most diagram of Fig. \ref{Prim-HZ8} has $\Cs=2$, the right-most $\Cs=8$, the others have $\Cs=1$.   
 

\subsection{The coefficients  $C_{n,[\alpha]}^{(g)} $}
We now return to our problem of determining the coefficients  $C_{n,[\alpha]}^{(g)} $ in (\ref{cumuleps}). 
From the previous discussion, if we denote $\CO_n([\alpha])\subset\CS_n$ the subset of  permutations of class $[\alpha]$, whose cycles involve only increasing sequences of integers, 
we have
\be\label{formul1} C_{n,[\alpha]}^{(g)} =
 \#\left\{\tau \big|\tau \in \CO_n([\alpha]), \ g =\oh\Big( n+1-\sum \alpha_\ell - \#\mathrm{cy}(\sigma\circ \tau)\Big)\right\}\,.\ee
Alternatively, one may use the diagrammatic language to write
\be\label{formul2} C_{n,[\alpha]}^{(g)} =\sum_{\mathrm{orbits}} \mathrm{length\ of\ orbit}=n   \sum_{\mathrm{orbits}} \frac{1}{\Cs} \,,\ee
\new{with a sum over orbits of diagrams of type $[\alpha]$ and genus $g$. }


\subsection{Remark on matrix integrals}
\label{matint}
As 't Hooft's double line notation~\cite{tH} suggests, 
the coefficient
\be\label{coeffC}  C_{n,[\alpha]}(\epsilon)=\sum_g C_{n,[\alpha]}^{(g)} \epsilon^g\ee  could be defined and computed  in  matrix integrals\\
-- (i) as  the coefficient  of $\prod_{\ell} \k_\ell^{\alpha_\ell} $ in the computation of $\langle \inv{N} :\tr M^n: \rangle_{rc}$
in a matrix theory with action $S= -\oh N \tr M^2+ N \sum_\ell  \new{\kappa_\ell}\, \tr M^\ell/\ell$; the notation $:\ :$ and the subscript ``{\it rc}" will be explained shortly; \\ 
-- (ii) as the value of  $\langle: \inv{N}  \tr M^n: \, :\prod_\ell\frac{(N \tr M^\ell/\ell)^{\alpha_\ell}}{\alpha_\ell!}:\rangle_{rc}$ in a Gaussian matrix theory.
\\
In both cases,  $\epsilon=\inv{N^2}$, if $N$ is the size of the (Hermitian) matrices;  $C_{n,[\alpha]}(N^{-2})$ is given by a sum of Feynman diagrams (in fact, of ``fat graphs",  or of maps) with $1+\sum_\ell \alpha_\ell$ vertices,
$n$ edges (``propagators") joining the $n$-vertex $\tr M^n$ to the other $\ell$-vertices,
and $f$ faces associated with each closed  index circuit. The double dots $:\ :$  is a standard notation in quantum field theory, 
where it denotes the {\it normal} or {\it Wick product}, that forbids edges from a
vertex to itself: here it forces all edges to reach the $n$-vertex. 
{The crucial point is that we impose a restricted crossing (``{\it rc}") condition: the edges connecting each $\ell$-vertex to the $n$-vertex cannot cross
one another, thus respecting their original cyclicity  and ordering. Only crossings of edges emanating from distinct vertices are allowed. 

 It is that constraint, a direct consequence of rule \ref{subsec21} (i) above, that makes the computation of the coefficients $C_{n,[\alpha]}^{(g)}$ 
by matrix integrals or group theoretical techniques, 
and the writing of recursion formulae between them, quite non trivial. For partitions into doublets, however, one deals only with 2-vertices for 
which the constraint is irrelevant, and $C_{n=2p, [2^p]}^{(g)}$ is computable by these techniques~\cite{WL1, HZ, 
Chapuy-th}.


 \begin{figure}[t]\begin{center}
{\includegraphics[width=.9\textwidth]{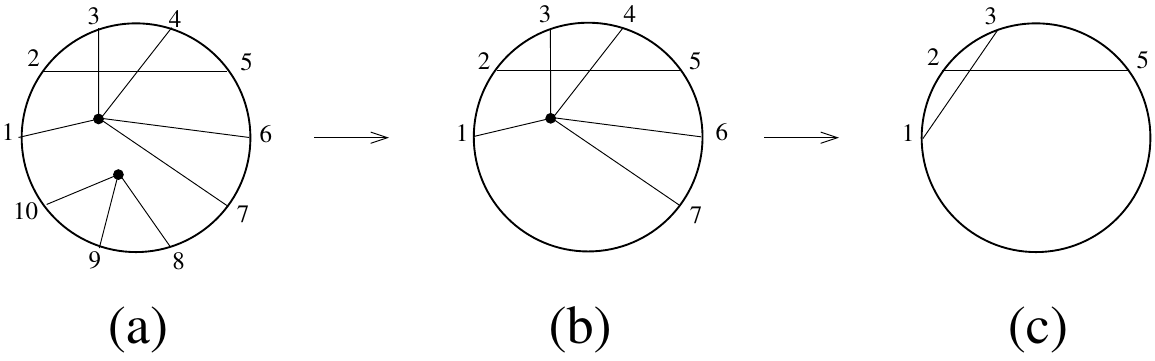}}
\end{center} 
\caption{\small{(a) Diagram for the partition of $\[10\]$  into $(\{1,3,4,6,7\},\{2,5\}, \{8,9,10\})$, $f=6$ hence genus $g=1-(3+1-10+ 6)/2=1$; (b) after 
removal of the three adjacent edges coming from the ``centipede" $ \{8,9,10\}$, here a
 3-vertex, now $n'=7$, $f'=4$, $g'=1$; (c) after reduction of two sets of adjacent edges to points 3 and 4, and 6, 7 and 1: now $n''=4$, $f''= 1$, $g''= 1+(2+1-4+1)/2=1$.   
 }}
\label{Fig2}
\end{figure}

 \begin{figure}[t]\begin{center}
{\includegraphics[width=.6\textwidth]{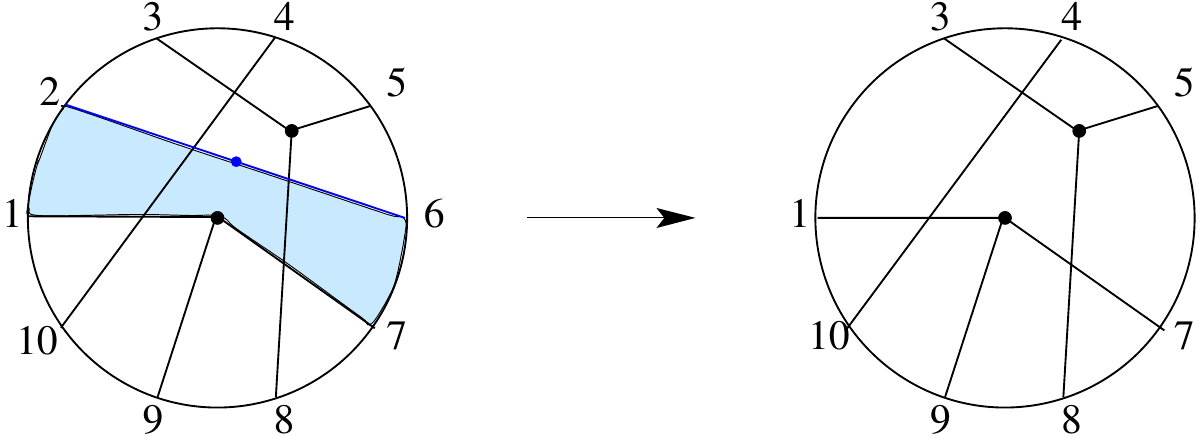}}
\end{center} 
\caption{\small{Removing the blue parallel pair of edges and the light blue face does not affect the genus: Variations $\Delta n=-2$, $\Delta f=-1$, $\Delta \sum \alpha_k=-1$,
hence $\Delta g=0$.   
 }}
\label{Fig3}
\end{figure}


\subsection{Reducing the diagrams}
\label{reduction}
In this subsection, we show that certain modifications of a diagram associated with a partition do not modify its genus.
The present discussion follows closely that of Cori and Hetyei~\cite{CoriH17}.
\\[4pt]
(i) Removing singletons.\\
Removing $p$ singletons changes the number of parts 
$\sum \alpha_k$ by $-p$, $n$ by $-p$ and the number of faces $f$ is unchanged, hence according to (\ref{Euler}) the genus remains unchanged.
 \\
{(ii) Removing centipedes.\\
{\bf Definition}. A {\it centipede} 
 is a planar linear subdiagram made of a $p$-vertex, all the  edges of which are attached in a consecutive way to the outer circle, 
Fig. \ref{Fig2}. In other words, it corresponds to a part of the partition  with consecutive integers (modulo $n$), $\{j, j+1,\cdots, j+p\}$. Removing it changes the number of parts 
$\sum \alpha_k$ by $-1$, $n$ by $-p$ and the number of faces $f$ by $-(p-1)$, see the figure, hence  the genus remains unchanged. \\}
(iii) Removing \new{adjacent} edges\\
If two  edges emanating from  a vertex  go to two consecutive points of the circle,  \new{(adjacent pair)}, 
see Fig \ref{Fig2}b-c, 
removing one of them does not change $\sum \alpha_k$  but changes $n$ and $f$ by $-1$, hence does not change the genus. One may 
iterate this operation on the same vertex until one meets a crossing with an edge emanating from another vertex. (If no crossing occurs, this means that the vertex \new{and its edges formed} a centipede  in the sense of (ii) and may be erased without changing the genus.) 
To allow an unambiguous reconstruction of all diagrams later in the dressing process, we adopt the following \\
{\bf Convention 1}: \new{in removing such adjacent edges, one keeps the edge attached to the marked point 1, or the first edge encountered
clockwise starting from 1, and one removes the others.}\\
\new{See Fig. \ref{Fig2} for illustration.}
 \\
 (iv) Removing parallel lines\\
{\bf Definition}. Two pairs of edges joining two vertices respectively to points $i$ and $j+1$, and  to points  $i+1$ and $j$ on the circle are
 said to be {\it parallel}. \\
 Note that this is equivalent to saying that they form a 2-cycle of the permutation $\sigma\circ \tau^{-1}$. And conversely, any such 
 2-cycle is associated with two parallel pairs of edges.
 \\ (a) If 
 one of these two vertices is a 2-vertex, one may remove the corresponding pair of edges \newa{and the 2-vertex} without changing the genus, since 
  $\sum \alpha_\ell$ and $f$  have decreased by 1 and $n$ by 2, see Fig. \ref{Fig3} for illustration. If both pairs of edges are attached to 
  2-vertices, we choose by {\bf Convention 2} to keep the pair attached to the point of the circle of  smallest label. In particular, if one of the pairs is
  attached to the marked point 1, it is kept and the other removed.\\
  (b) If both pairs of edges are attached to vertices of valency larger that 2, we keep them both. See Fig. \ref{SemiPrim10-3322} below for an example.\\
  After carrying these removals of parallel lines, we are left with {\it primitive} or {\it semi-primitive} diagrams (or partitions),
  following Cori--Hetyei's terminology: in primitive diagrams, no  parallel pair is left; therefore, by the remark above, all
  cycles of  $\sigma\circ \tau^{-1}$ have length larger than 2. Semi-primitive diagrams  still have parallel pairs of type (b).

\bigskip
Now Cori and Hetyei have proved some fundamental results:\\
{\bf Proposition}. {\sl To an arbitrary diagram corresponds a unique primitive (or semi-primitive) diagram obtained by a sequence of reductions as
above, and independent of the order of these reductions.}\\
Our new observation is that,  conversely, any diagram may be recovered by ``dressing"  a primitive (or semi-primitive) diagram, as we shall
see below.

Moreover, \\ 
{\bf Proposition}. \cite{CoriH17} {\sl For a given genus, there are only a finite number of primitive diagrams.} \\
This follows from  two inequalities:  $f=\#\mathrm{cy}(\sigma\circ\tau^{-1}) \le\frac{ n}{3}$, since in a primitive diagram all cycles of $\sigma\circ\tau^{-1}$
are of length larger or equal to three (see above); and $\sum \alpha_i \le n/2$ after eliminating the singletons. Hence plugging these inequalities in 
(\ref{genus}), we get for a primitive diagram
\be\label{nprim}  n\le 6(2g-1)\,. \ee
As for the semi-primitive diagrams, it was shown in \cite{CoriH17} that they are all obtained by a finite number of operations 
from the primitive ones, hence are themselves in finite number.

For example\\
{\bf Proposition}\cite{CoriH17} {\sl The irreducible diagrams  of genus 1 are the two diagrams of Fig. \ref{Fig4}, that have two, resp. three 2-lines. 
\new{No semi-primitive occurs in genus 1.}}
\\
The proof of that  statement is given in \cite{CoriH17}, sect. 8, where the two \new{primitive} partitions or diagrams are referred to as $\beta_1$ and $\beta_2$. \\
\moi{\tiny Checked in Graphic2.nb}


\section{From genus 0 to genus 1 \dots}

 \begin{figure}[t]\begin{center}
{\includegraphics[width=.15\textwidth]{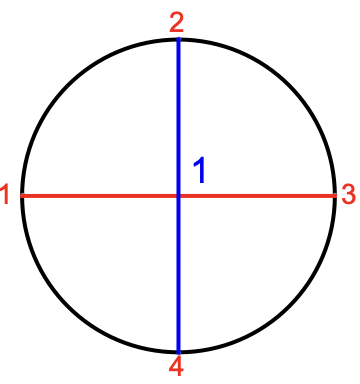}\qquad\qquad \raisebox{2pt}
{\includegraphics[width=.15\textwidth]{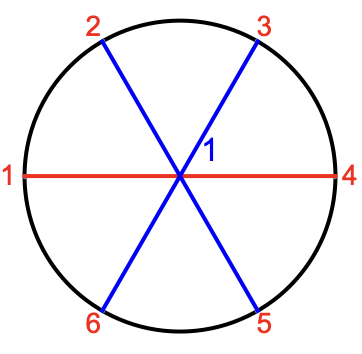}}}
\end{center} 
\caption{\small{The two ``primitive" diagrams of genus 1. The blue figure in the middle is the weight of the diagram in (\ref{formul1}), 
namely the length of its orbit.
 }}
\label{Fig4}
\end{figure}

\subsection{Non crossing partitions and the genus 0 generating function}
Recall first that in genus 0, the formula given by Kreweras \cite{Krew} on the census of non crossing 
partitions may be conveniently encoded in the following functional relation between the genus 0
GF of moments $Z^{(0)}(x)$ and that of cumulants $W(x)$ defined above \footnote{  Recall this relation is equivalent to the functional identity 
$P\circ G=\mathrm{id}$, where $G(u):= u^{-1} Z^{(0)}(u^{-1})$ and $P(z):= z^{-1} W(z)$,
and $R(z)=P(z)-\inv{z}$ is the celebrated Voiculescu $R$ function~\cite{V86, Speicher}. }
\be\label{Z0W} Z^{(0)}(x) = 1+ W(x Z^{(0)}(x) )\,. \ee
Indeed by application of Lagrange formula, one recovers Kreweras' result
\be \label{genus00}   
C^{(0)}_{n,[\alpha]}= \new{\frac{n!}{(n+1-\sum \alpha_k)!  \ \prod_k \alpha_k!}}\,, 
 \ee
as proved in \cite{BIPZ}. 

There is a simple diagrammatical interpretation of the relation (\ref{Z0W})
due to Cvitanovic \cite{Cvitanovic}, 
see Fig.~\ref{cvitanovic1},
which reads: {\it in \newa{an arbitrary} planar (i.e., non-crossing) diagram, the marked point 1 on the exterior circle is necessarily
connected to a $n$-vertex, $n\ge 1$, between the $n$ edges of which 
lie arbitrary insertions of other \new{(linear)} diagrams of $Z^{(0)}$}. \\ Our aim is  to extend this 
kind of relation to higher genus.

 \begin{figure}\begin{center}
{\includegraphics[width=.7\textwidth]{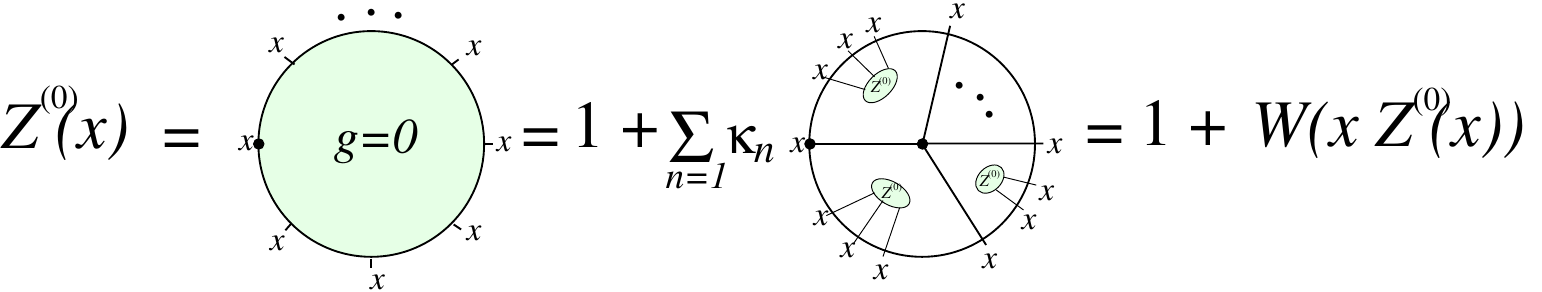}}
\end{center} 
\caption{\small{ A graphical representation of identity $Z^{(0)}(j)= W(j\, Z^{(0)}(j) )$  }}
\label{cvitanovic1}
\end{figure}

\begin{figure}
\begin{center}
{\includegraphics[width=.95\textwidth]{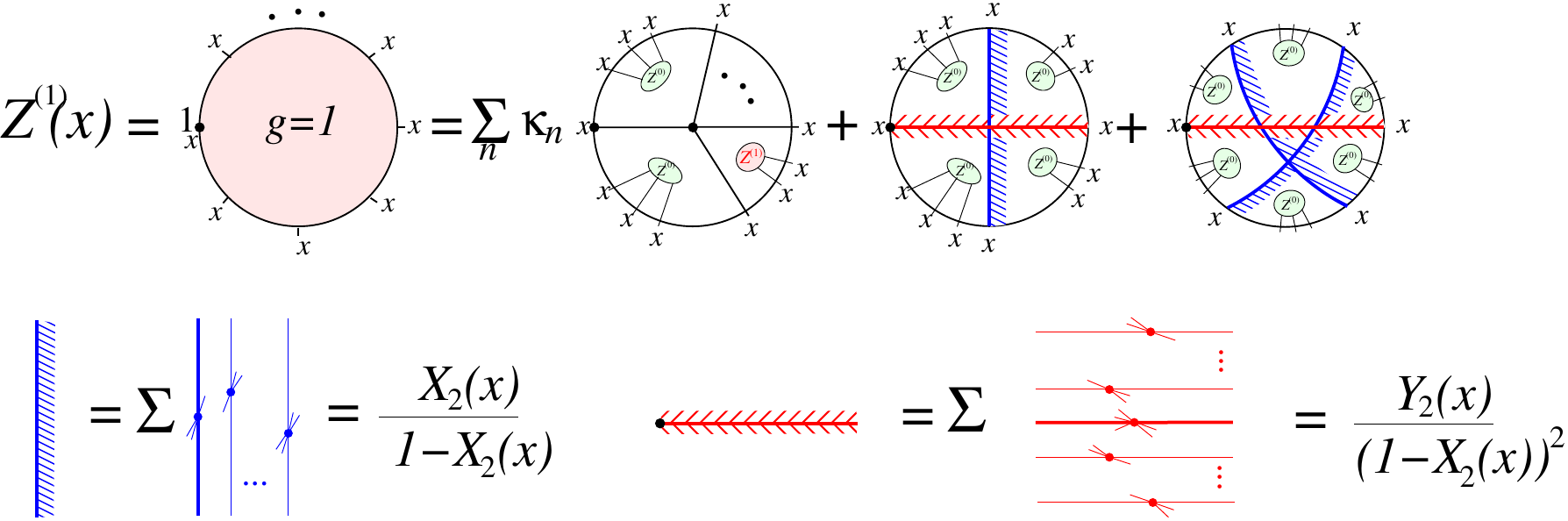}}
\end{center} 
\caption{\small{Top: A graphical representation of identity (\ref{identg=1}).  
Bottom: a schematic representation of the dressing of the (red) 2-line attached to the marked point 1; or to another
(blue) 2-line. In the latter case, according to  Convention 1,  additional edges may be ``emitted" from the central vertex 
to go to clockwise adjacent points on the circle, and their contribution to the generating function is $X_2(x)$. 
 For the red line, these additional  edges may connect 
to either side of the marked point, and they contribute  $Y_2(x)$ to the GF. 
}}
\label{genus1}
\end{figure}


\subsection{Dressing the genus 1 primitive diagrams}
\label{dressingofg=1}

We have seen that genus 1 diagrams may be reduced to the two primitive ones of Fig. \ref{Fig4}. 
 We  now write  a  relation \`a la Cvitanovic between the generating functions $W$, $Z^{(0)}$ and $Z^{(1)}$, depicted in 
Fig.~\ref{genus1}
\be\label{relg=1} Z^{(1)}(x)= \sum_{n\ge 2} \kappa_n n x^n (Z^{(0)} )^{n-1} Z^{(1)} + 
\mathrm{sum\ of\ dressed\ diagrams\ of\ Fig.\ \ref{Fig4}}\,,
\ee 
which reads: in a generic diagram of genus 1, 
the marked point 1 is attached {\bf (a)}  either to   an edge of an $n$-vertex, between the non-crossing edges
of which are inserted one (linear) subdiagram of genus 1 and $(n-1)$ subdiagrams of genus 0~\footnote{Remember that by convention,
$Z^{(0)}(x)$ starts with 1, hence these subdiagrams of genus 0 may be trivial}, {\bf (b)}  or to an edge of a 
dressed primitive diagram of genus 1.

Let us concentrate on the case {\bf (b)} and make explicit what is meant by dressing.

The dressing consists in reinserting the elements removed in steps (iv)--(i) of sect. \ref{reduction}, in that reverse order.
First, additional 2-lines  are introduced, ``parallel" to the two, resp. three 2-lines  
of the primitive diagrams of Fig.~\ref{Fig4}. Each of these 2-lines carries by definition a 2-vertex.  Then to reinsert
``adjacent" edges removed in step (iii), each of these 2-vertices may be transformed into a $k$-vertex, whose $k-2$ additional edges may fall,
by Convention 1, on either of the two arcs of the circle adjacent to the end points of the 2-line
 and ``clockwise downstream", and without crossing one another: 
 there are  $k-1$ partitions of $k-2$ into two parts, one of them possibly empty, hence we attach a weight $X_2(x):= \sum_{k\ge 1}(k-1)  \kappa_k x^k $
 to each of these parallel lines. 
Since there is an arbitrary number $r\ge 0$ of  parallel lines, they contribute  $X_2(x)^r$, and their geometric series sums up to $1/(1-X_2(x))$.
The same applies to the  original blue 2-lines of the primitive diagram of Fig.~\ref{genus1}, which thus gives each a factor $X_2(x)$. The red 2-line, which is the one attached to the marked point 1, has a different weight, as 
the $k-2$ edges emanating from its $k$-vertex may fall  on either side of the marked point or on the rightmost part of the 
diagram (see Convention 2 above): this is associated with a partition of the $k-2$ edges into three parts (two of them possibly empty), 
in number $k(k-1)/2$,
which gives the red  2-line a weight $Y_2(x)=\sum_k  \frac{k(k-1)}{2}  \kappa_k x^k$, while its dressing by  parallel lines leads to a factor $1/(1-X_2(x))^2$, because again, parallel lines above or below the red 2-line are possible. Last step   
consists in reinserting ``centipedes" and (possibly) singletons, 
namely in changing everywhere $x$ into $\tilde x=x Z^{(0)}(x)$. \\ 
 In that way, we have reinstated all features that had been erased in the reduction to primitive diagrams,   and constructed the 
 contribution to the GF $Z^{(1)}(x)$ of {\it all diagrams in which the marked point 1 is attached to an edge that belongs to a 
 dressed primitive diagram.} 
 Indeed in the resulting diagrams, the marked point 1 may be attached to any of the edges, as it should: this is
 clear whenever that edge is an edge of the primitive diagram; this is also true if the edge is  one of the parallel lines added, or  one of the added 
 adjacent edges: that was the role of the factors in the definition of $X_2$ or $Y_2$ to count these cases. 
   It is thus clear that all possible diagrams of type {\bf (b)} contributing to  $Z^{(1)}$ have been obtained by the dressing procedure, and that they are generated once and only once, hence with the right weight.
   Finally the cases {\bf (a)} where 1 is not attached to a dressed primitive,
 but to some genus 0 subdiagram,  are accounted for  by the first  term  in equ.(\ref{relg=1}). \\


\subsection{The genus 1 generating function}
\label{GFg=1}
\omit{We may now write  a Cvitanovic-like relation between the generating functions $W$, $Z^{(0)}$ and $Z^{(1)}$, depicted in Fig.\ref{genus1}
\be\label{relg=1} Z^{(1)}(x)= \sum_{n\ge 2} \kappa_n n x^n (Z^{(0)} )^{n-1} Z^{(1)} + 
\mathrm{sum\ of\ dressed\ diagrams\ of\ Fig.\ \ref{Fig4}}\,,
\ee 
which reads: the marked point 1 is either attached to an edge of a $n$-vertex, between the non-crossing edges
of which are inserted one (linear) subdiagram of genus 1 and $(n-1)$ subdiagrams of genus 0~\footnote{Remember that by convention,
$Z^{(0)}$ starts with 1, hence these subdiagrams of genus 0 may be trivial}, or to an edge of a 
dressed primitive diagram of genus 1.  }

Define  
\be\label{deftx} \tx=x Z^{(0)}(x)\,.\ee
Gathering all the contributions of sect.~\ref{dressingofg=1} we have 
\be\label{identg=1} Z^{(1)}(x)= \sum_{n\ge 2} \kappa_n n x^n (Z^{(0)}(x) )^{n-1} Z^{(1)}(x) + \frac{Y_2(\tilde x)X_2(\tilde x)}{(1-X_2(\tilde x))^3}
+ \frac{Y_2(\tilde x)X_2^2(\tilde x)}{(1-X_2(\tilde x))^4}\,, \ee
\ie $$ (1-V(x)) Z^{(1)}(x)=  \frac{Y_2(\tilde x)X_2(\tilde x)}{(1-X_2(\tilde x))^4}$$
where \bea \label{defX2} X_2(x)&=&  \sum_{k\ge 2}(k-1) \kappa_k x^k = x W'(x)-W(x)\,, \\
\label{defY2}  \X1(x)&=& \sum_{k\ge 2}\frac{k(k-1)}{2}\kappa_k x^k= \oh x^2 W''(x)\\    
 \label{defV} V(x)&=& \sum_k k \kappa_k x^{k} Z^{(0)\, k-1}   =x W'(\tilde x)\,. \eea 
This is summarized in the following theorem.

\bigskip
{\bf Theorem 1.}\label{Thm1} {\sl  If  $\tx=x Z^{(0)}(x)$, the generating function of genus 1 partitions is given by }
\be\label{Zgenus1} 
Z^{(1)}(x) = \frac{ X_2(\tilde x)  \X1(\tilde x)}{(1-X_2(\tilde x))^4\,(1-V(x))} \,.\ee
Alternatively, if we introduce
\be \label{deftX2tY2}   \tX_2(x):=\frac{X_2(\tx)}{(1-X_2(\tx))}\qquad {\tY}_2(x):=\frac{\X1(\tx)}{(1-X_2(\tx))^2}\ee
we have the simple expression
\be\label{Z1alt}  Z^{(1)}(x)=  \frac{ \tY_2(x) \tX_2(x)(1+\tX_2(x))}{(1-V(x))} \ee

 
\subsection{Examples and applications}
\subsubsection{$n=2p$, $ [2^p]\vdash n$}
If all $\kappa_i$ vanish but $\kappa_2=1$, \ie if we consider partitions of $n=2p$ into $p$ doublets, which is the celebrated case 
considered in \cite{WL1, HZ}, we have $W(x)=x^2$, hence 
\be\label{HZg0} Z^{(0)}(x; \kappa_2=1, \kappa_{i\ne 2}= 0)=   \frac{1-\sqrt{1-4x^2}}{2x^2} \ee
as the solution of equ. (\ref{Z0W}).  Then following Theorem 1,   we find
\be\label{HZg1}  Z^{(1)}(x; \kappa_2=1, \kappa_{i\ne 2}= 0)= \frac{x^4}{(1 - 4 x^2)^{5/2}}\,,\ee
in accordance with known results.
  
\subsubsection{$n=3p$, $[3^p]\vdash n$}
\label{nto3p}
In that case, we take $\kappa_3=1$, $W(x)=x^3$, hence $Z^{(0)}$ satisfies the third degree equation, 
\be\label{cubicequ} (x Z)^3-Z+1=0\ee 
and it is the GF of Fuss--Catalan numbers. We may write it as 
\be\label{FCg0} Z^{(0)}(x; \kappa_3=1, \kappa_{i\ne 3}= 0)  =\frac{2}{\sqrt{3 x^3}} \sin\Big(\frac{1}{3}\Arcsin\big(\frac{3}{2}\sqrt{3x^3}\big)\Big)\,.%
\ee
Then following Theorem 1,  one finds, after some algebra,
\be\label{FCg1}Z^{(1)}(x; \kappa_3=1, \kappa_{i\ne 3}= 0)  =\frac{1152\, x^3 \sin ^6\left(\frac{1}{3} \Arcsin\big(\frac{3 \sqrt{3x^3}}{2}\big)\right)}
{\left(2 \cos \left(\frac{1}{3} \Arccos\big(1-\frac{27 x^3}{2}\big)\right)-1\right)
\left(9 \sqrt{x^3}-4 \sqrt{3} \sin   \left(\frac{1}{3} \Arcsin\big(\frac{3 \sqrt{3x^3}}{2}\big)\right)\right)^4 }\ee
with a Taylor expansion
$$6 x^6 + 102 x^9 + 1212 x^{12} + 12330 x^{15} + 114888 x^{18} + 1011486 x^{21} + 
 8558712 x^{24} + 70324884 x^{27} + 564931230 x^{30}+\cdots$$
 \new{in agreement with direct calculation,} see  
 \cite{comp}. Note that the closest singularity of $Z^{(1)}$ is at the vanishing point of the discriminant of (\ref{cubicequ}),
 namely $x^3=4/27$: 
 $$ Z^{(1)}(x; \kappa_3=1, \kappa_{i\ne 3}= 0) \sim \frac{\mathrm{const.}}{(\frac{4}{27}-x^3)^{5/2}}\,,$$
 \new{when $x^3\to 4/27$, with the same exponent $5/2$ as in (\ref{HZg1})}.

 \subsubsection{Total number of partitions of genus 0 and 1}
 \label{totalnbr01}
Let all $\kappa$ be equal to 1, resp. all $\kappa$'s but $\kappa_1=0$. Then the previous expressions yield the GF of the numbers of partitions of genus 0 or 1, 
with, resp. without singletons:
\bea
\label{Z0} 
Z^{(0)}(x; \kappa_i=1) &=&    \frac{1-\sqrt{1-4x}}{2x} \\ 
\nonumber 
 \hat Z^{(0)}(x):= Z^{(0)}(x; \kappa_1=0, \kappa_{i\ge 2}=1) &=&    \frac{1-\sqrt{1-\frac{4x}{1+x}}}{2x} =\frac{1+x-\sqrt{1-2x -3x^2}}{2x(1+x)}\qquad \mathrm{no\ singleton}\\
\label{Z1}  
Z^{(1)}(x; \kappa_i=1) &=&  \frac{x^4}{(1-4x)^{5/2} }   \\ 
\label{Z1p}  
 \hat Z^{(1)}(x):= Z^{(1)}(x; \kappa_1=0, \kappa_{i\ge 2}=1)  &=&    \frac{x^4}{(1-2x-3x^2)^{5/2}} \qquad  \mathrm{no\ singleton}
\eea
on which we may verify the relations (\ref{ZZhat}-\ref{ZhatZ}) above.
\omit{\blue [and that we may supplement for genus 2 with the result of Cori and Hetyei \cite{CoriH17}
$$ Z^{(2)}(x) = \frac{x^6 (1 + 6 x - 19 x^2 + 21 x^3)}{(1 - 4 x)^{11/2}}\,. $$
and   $$\hat Z^{(2)}(x)=\frac{x^6 (1 + 10 x + 5 x^2 + 5 x^3 + 9 x^4)}{(1 - 2 x - 3 x^2)^{11/2}}$$ as a consequence of (\ref{ZZhat}).
\moi{\tiny Checked in FromGenus0ToGenus1ToGenus2.nb} ]}

{\it Proof}. If all $\kappa_i=1$, $W(x)=x/(1-x)$ as a formal series, and $Z^{(0)}(x)$, solution of $ Z^{(0)}(x)= W(x Z^{(0)}(x))$ as a formal series, is
given by (\ref{Z0}), (the GF of the Catalan numbers). Likewise, if $\kappa_1=0$, the others equal to 1, $W(x)=x^2/(1-x)$, etc. For genus 1, we then 
make use of Theorem 1 to derive (\ref{Z1}-\ref{Z1p}). \qed

\subsubsection{Number of partitions with a fixed number of parts, in genus 0 and 1}
\label{FGxy01}
Let all $\kappa$ be equal to $y$, then $W(x)= xy/(1-x)$, and $Z^{(g)}(x,y)= \sum_{n,k}  p^{(g)}(n,k)   x^n y^k$ 
is the GF of the numbers $p^{(g)}(n,k)$ of genus $g$ partitions of $n$ with $k$ parts. $Z^{(0)}$ is the solution of equ. (\ref{Z0W})
\be  \label{fixednbrparts0} Z^{(0)}(x, y)= \frac{1 + x- x y - \sqrt{ (1 + x - x y)^2-4x}}{2 x}\,. \ee
which is the GF of Narayana numbers, and then we compute by (\ref{Zgenus1})
\be  \label{fixednbrparts1} Z^{(1)}(x,y)=  \frac{x^4 y^2}{( (1 +x -x y)^2 -4x)^{5/2}}\ee
which is the expression given by Yip~\cite{Yip}, and Cori and Hetyei \cite{CoriH13}.

If we exclude singletons, $W(x; \kappa_1=0)=x^2y/(1-x)$, and the GF read now
\bea  \label{fixednbrpartsnosingle} \hat Z^{(0)}(x,y):= Z^{(0)}(x, y; \kappa_1=0)&=&\frac{{1 + x}- \sqrt{(1 -  x)^2  - 4 x^2 y}}{2 x(1+x y )} \\
\nonumber  \hat Z^{(1)}(x,y):=Z^{(1)}(x, y; \kappa_1=0)&=&  \frac{x^4 y^2}{((1 -  x)^2  - 4 x^2 y)^{5/2}}\,.
\eea


\section{\dots to genus 2}

\subsection{Primitive and semi-primitive diagrams of genus 2}
The list of primitive and semi-primitive diagrams of genus 2 is known, thanks to the work of Cori and Hetyei \cite{CoriH17}. 
This has been confirmed independently, in the present work, by generating on the computer all partitions of genus 2 of a given type, and
then eliminating all those that involve adjacent or parallel edges. By inequality (\ref{nprim}) these primitive diagrams
have at most 18 points (\ie $n\le 18$), and either up to 9 2-vertices, or one or two 3-vertices, or one 4-vertex.
In Table 1, are listed their number for increasing total number of points $n$.~\footnote{In Table 1 of \cite{CoriH17} there is the unfortunate omission 
of the 175 primitive diagrams with one 3-vertex (a 3-cycle in their terminology), while those diagrams are properly taken into account 
in the ensuing formulae. These missing diagrams are listed in Fig. \ref{Prim15-32222222}.}\\

\begin{center}\begin{tabular}{l | c| c|  c| c|  c| c|  c| c|   }
\qquad & 2-vertices & one 3-vertex & two 3-vertices & two 3-vertices & one  4-vertex  \\
$n$  &  & & & semi-prim.&\\
\hline
6 & 0 & 0  & 1 & 0 & 0  \\
7 & 0 & 14& 0 & 0 & 0\\
 8 &21 & 0 & 20& 0& 6 \\
 9 &  0& 141& 0& 0&0\\
10 &168& 0& 65& 15 &15 \\
11  &0& 407& 0& 0&0 \\
12 &483& 0& 52&36& 9 \\
13  &0& 455& 0& 0&0 \\
14  &651& 0& 0& 21&0 \\
15  &0& 175& 0& 0&0\\
16  &420& 0& 0& 0&0\\
17  &0& 0& 0& 0&0 \\
18  &105& 0& 0& 0&0\\
\hline
\end{tabular}\\[6pt]\end{center}
\centerline{Table 1. Number of (semi-)primitive diagrams of genus 2.}

\bigskip
Based on this list of primitive diagrams, we may now write an equation similar to (\ref{identg=1})

\bea\label{Zgenus2}  Z^{(2)}(x)&=& \sum_n  n \kappa_n x^n (Z^{(0)}(x))^{ n-1}  Z^{(2)}(x)  
\\  &+& \nonumber
\quad\mathrm{dressing\ of\ (semi-)primitive\ diagrams\ of\ genus\ 2} 
 \eea
 as illustrated in Fig. \ref{genus2}.
\\ \newa{{\bf Remark.} It might seem natural to also have in the r.h.s. of (\ref{Zgenus2} ) a term with two 
insertions of genus 1 subdiagrams. In fact such diagrams will be included in the set  of primitives and their
dressings. An example is given by the first diagram of Fig.~\ref{Prim-HZ8}. }

  \begin{figure}\begin{center}
{\includegraphics[width=.85\textwidth]{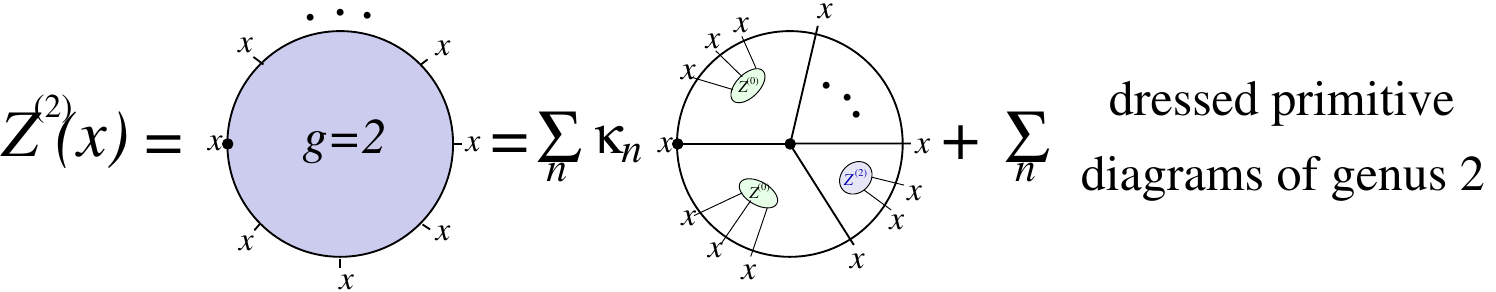}}
\end{center} 
\caption{\small{ A graphical representation of relation (\ref{Zgenus2}) } } 
\label{genus2}
\end{figure}

\def\kk2{\kappa_2}\def\kk2{}


\subsection{Dressing  of primitive diagrams of genus 2}
The dressing of primitive diagrams with only 2-lines (Column 2 of Table 1) involves 
the same functions $\tX_2$ and $\tY_2$ defined above in sect. \ref{GFg=1}: 
$\tY_2$ is assigned to the line attached to point 1, while  the other lines carry the weight $\tX_2$. 
Hence  their contribution to the r.h.s. of  equ.(\ref{Zgenus2}) reads
 \be\nonumber z_2=   \tY_2(x) \Big(21 \tX_2^3(x) + 168 \tX_2^4(x) 
  +483 \tX_2^5(x) + 651   \tX_2^6(x) + 420  \tX_2^7(x) +   105  \tX_2^8(x) \Big)\ee
  with the notations of (\ref{deftX2tY2}). 

 For the dressing of primitive diagrams with 3- or 4-vertices, we must introduce new functions
that generalize $X_2$ and $Y_2$ defined in (\ref{defX2}-\ref{defY2})
 \bea\label{newfns} 
 X_\ell(x) &=& \sum_{k\ge \ell}  {k-1\choose \ell-1} \kappa_k x^k\\
 \nonumber Y_\ell(x) &=& \sum_{k\ge \ell}  {k\choose \ell} \kappa_k x^k\\
  \nonumber\ell>2 \qquad \tX_\ell(x) := \frac{X_\ell(\tx)}{(1-X_2(\tx))^{\ell}}\quad &;&\quad  \tY_\ell(x) := \frac{Y_\ell(\tx)}{(1-X_2(\tx))^{\ell}}\,.
 \eea
with, as before, $\tx = x Z^{(0)}(x)$. (Beware that the power of $(1-X_2(\tx))$ in the denominator of $\tX_\ell$ does not apply to $\ell=2$, compare with (\ref{deftX2tY2}).)  These functions too may also be expressed in terms of derivatives of $W$:  for example, $Y_3(x)=\inv{6}x^3 W'''(x)$, etc.

Consider first a primitive diagram with one 3-vertex, like those depicted in Fig. \ref{Prim7-322}. Remember that all distinct rotated diagrams
must be considered and hence, the marked point 1 may be attached  to the 3-vertex or to any one of the  2-lines. \\
(i) In the case where the marked point 1 is attached to one of the 2-lines, its 2-vertex may be changed into a $k$ vertex, $k> 2$ and
as in sect.  \ref{dressingofg=1}, this yields a weight $Y_2(x)/(1-X_2(x))^2$, while the lines emanating from the 3-vertex or parallel to it
contribute $X_3(x)/(1-X_2(x))^3$. And again, a final change of $x$ into $\tx$ completes the dressing.\\
(ii) In the former case, 1 attached to the 3-vertex, 
this 3-vertex may be promoted to a $k$-vertex, $k>3$, with $k-3$ lines ending on four different arcs of the circle: there are ${k\choose 3}$ ways 
of distributing them,  whence a weight $Y_3(x)$. Then adding parallel lines may be done in 3 ways, whence a weight $1/(1-X_2(x))^3$.
The  2-lines, on the other hand, carry a weight $X_2(x)/(1-X_2(x))$, just like   in sect.  \ref{dressingofg=1}.
Finally, again as in sect. \ref{dressingofg=1}, the variable $x$ has to be substituted for the dressed one $\tx= x Z^{(0)}$ to take into 
account all possible insertions of genus 0 subdiagrams. \\
(iii) There is, however, a case not yet accounted for by the previous dressing.  When the marked point 1 is attached to a 2-line parallel 
to a pair of edges of the 3-vertex, that line has been erased in the reduction process and must be restored. A weight $2Y_2(x)/(1-X_2(x))$ is
attached to that new line, with a factor 2 comes from the two ends of the 2-line, and  a single factor $1/(1-X_2(x))$ as compared with what we saw
in sect. \ref{dressingofg=1}, because the counting of parallel lines between the new line and the 3-vertex has already been taken into account  in the term $\tX_3(x)$.  

Now each of the previous contributions must be weighted by its number of occurrences when the diagram is rotated.
For example, each of the two diagrams of Fig. \ref{Prim7-322} contributes 
 $+\ 4 \tY_2 \tX_2\tX_3$ (since marked point 1 may be at any of the four end-points of the 2-lines) +$\  3 \tY_3 \tX_2^2$ (3 ways of attaching point 1 to the 3-vertex)
$\ + 3 \tX_3 \tX_2^2(2Y_2(\tx)/(1-X_2(\tx))  $  (when 1 lies on a line parallel to two edges of the 3-vertex). More generally, for a 
primitive diagram of an orbit of symmetry order $\Cs$,  with one 3-vertex and $p$ 2-lines, $n=3+2p$, the weight is
$$\inv{\Cs} \left(2p  \tY_2 \tX_2^{p-1} \tX_3 + 3 \tY_3 \tX_2^{p}  + 3 \tX_3 \tX_2^{p} (2Y_2(\tx)/(1-X_2(\tx))\right)  \,,$$
where we write $\tX_\ell$ and $\tY_\ell$ in short for $\tX_\ell(x)$ and $\tY_\ell(x)$.
Thus the orbits of partitions of $\[n\]$ with a primitive diagram with a single 3-vertex
contribute $$ \sum_{\mathrm{orbits}}\inv{\Cs} \left((n-3) \tY_2 \tX_2^\frac{n-5}{2} \tX_3 
+ 3  \tX_2^{\frac{n-3}{2} } \Big(\tY_3 + \tX_3  \frac{2Y_2(\tx)}{(1-X_2(\tx))}\Big)\right)  \,. $$
But as we saw in (\ref{formul2}),  for a given $n$, $\sum_{\mathrm{orbits}}\inv{\Cs} =\frac{N}{n}$, where $N$ is the number listed in Table 1, column 3, row $n$.
In total the diagrams with a single 3-vertex contribute to the r.h.s. of (\ref{Zgenus2}) 
 the amount $z_3$ listed below in (\ref{z3}).

The dressing of primitive diagrams with two 3-vertices  or one 4-vertex (columns 4 and 6 of Table 1) is done along 
similar lines. Thus for an orbit of  primitive diagram with two 3-vertices and $p$ 2-lines, with now $n=2p + 6$,  we get
\be\label{2-3v}\inv{\Cs} \left(2 p \tY_2   \tX_3^2   \tX_2^{p-1}  
+6 \tX_2^{p } \tX_3 \Big(\tY_3 + \tX_3  \frac{2Y_2(\tx)}{(1-X_2(\tx))}\Big)\right) 
\ee
and the total contribution $z_{33}$ of such diagrams is given in  (\ref{z33}). \\

For a primitive diagram with one 4-vertex and $p$ 2-lines, (and $n=2p+4$), likewise, we get 
$$ \inv{\Cs} \left( 2p  \tY_2   \tX_4   \tX_2^{p-1}   
+4 \tX_2^{p }\Big(\tY_4 + \tX_4  \frac{2Y_2(\tx)}{(1-X_2(\tx))}\Big)\right) $$
and the total contribution $z_4$ is given in  (\ref{z4}). \\

Finally, the dressing of semi-primitive diagrams (see a sample in Fig. \ref{SemiPrim10-3322}) requires special care to avoid double counting.  
Consider such a semi-primitive diagram, thus  with two 3-vertices and $p$ 2-lines, $n=2p+6$. 
First, when the point 1 is attached to one of the 2-lines or one of the two
3-vertices, we have a contribution  like the first two terms in (\ref{2-3v}), but multiplied by $(1-X_2(\tx))$ not to count 
twice the set of lines between the two parallel lines. Moreover, when the point 1 is attached to  an added
line parallel to one of the branches of the two 3-vertices, there are 5 locations for that line, whence a contribution
$\frac{5}{\Cs} \tX_3^2 \tX_2^p \times 2Y_2(\tx)$, with no further factor $1/(1-X_2(\tx))$. In total, a semi-primitive diagram contributes
$$ \inv{\Cs} \left((1-X_2(\tx)) \Big(2 p \tY_2   \tX_3^2   \tX_2^{p-1} 
+ 6 \tY_3 \tX_3  \tX_2^p \Big)
+5  \tX_2^{p } \tX_3^2 (2Y_2(\tx))\right) $$
 and the total from semi-primitive diagrams appears as  $z_{33s}$ in (\ref{z33s}). \\[4pt]

\neww{{\bf Remark}. As noticed by Cori and Hetyei~\cite{CoriH17}, the semi-primitive diagrams may be obtained from the primitive ones by 
``splitting" a vertex of valency larger than 3. For example the three diagrams of Fig.~\ref{SemiPrim10-3322} may be obtained from those of Fig.~\ref{Prim8-422} 
by splitting their 4-vertex as in Fig.~\ref{splitting}. One might thus consider only primitive diagrams and include the splitting operation in the 
dressing procedure. The benefit is that primitive diagrams are easy to characterize: they are such that, {\it in genus 2},  the permutation $\tau$ has no 
1-cycle and $\sigma\circ\tau^{-1}$ no 2-cycle.  }


\subsection{General case of genus 2}
Collecting all the contributions of the previous subsection, we can now make equation (\ref{Zgenus2}) more explicit in the form of \\[10pt]
{\bf Theorem 2.}\label{Thm2} {\sl The generating function of genus 2 partitions is given by }
\be \label{Z2genus2} \qquad \qquad Z^{(2)}(x)(1-V(x))= z_2+ z_3+ z_{33}+z_{33s} + z_4 \ee
{\sl where $V(x)$ has been given in (\ref{defV}) and $z_2,\cdots, z_4$ are the contributions of dressing the (semi-)\-primitive diagrams listed in Table 1.}
\bea 
 \label{z2}  z_2&=&   \tY_2 (21 \tX_2^3 + 168 \tX_2^4 + 483 \tX_2^5 + 651   \tX_2^6 + 420  \tX_2^7   + 105  \tX_2^8 )\,;\\[4pt]
 \label{z3} z_3&=& \tX_3  \tY_2 (8 \tX_2 + 94 \tX_2^2 + 296 \tX_2^3 + 350 \tX_2^4 + 140 \tX_2^5) \\
      \nonumber&&\qquad  +\tX_2 (6 \tX_2  + 47 \tX_2^2 + 111 \tX_2^3  + 105 \tX_2^4 +35 \tX_2^5)\Big(\tY_3 + \tX_3  \frac{2Y_2(\tx)}{(1-X_2(\tx))}\Big)\,; \\
 \label{z33} z_{33} &=& \tX_3^2 \tY_2 (5 + 26 \tX_2 + 26 \tX_2^2   ) \\
     \nonumber  && \qquad +  \tX_3 (1 + 15 \tX_2 + 39 \tX_2^2 + 26 \tX_2^3) \Big(\tY_3 + \tX_3  \frac{2Y_2(\tx)}{(1-X_2(\tx))}\Big)\,; \\
 \label{z33s}   z_{33s} &=& \tY_2  \tX_3^2 \tX_2 (6 + 18 \tX_2 + 12 \tX_2^2 )(1 - X_2(\tx))  \\ 
 \nonumber   &&\qquad +  \tY_3  \tX_3  \tX_2^2 (9 + 18 \tX_2 + 9 \tX_2^2) (1 - (X_2(\tx)) 
   +\tX_3^2 \tX_2^2 (15 + 30 \tX_2 + 15 \tX_2^2)  Y_2(\tx)\, ;\\
   \label{z4}   z_4& =&  \tY_2 \tX_4  (3  \tX_2 + 9  \tX_2^2 + 6 \tX_2^3) 
   +    (3 \tX_2^2 + 6 \tX_2^3 + 3 \tX_2^4) \Big(\tY_4 + \tX_4  \frac{2Y_2(\tx)}{(1-X_2(\tx))}\Big)\,,
 \eea
and we recall that $\tX_\ell$ and $\tY_\ell$ stand for  $\tX_\ell(x)$ and $\tY_\ell(x)$ defined in (\ref{newfns}).

The resulting expressions for the numbers $C_{n,[\alpha]}^{(2)}$ have been tested up to $n=15$ and all $[\alpha]$ against
  direct enumeration using formulae (\ref{formul1}) or (\ref{formul2}), and for some higher values of $n$ for a few particular cases.

    \begin{figure}
\begin{center}
{\includegraphics[width=.85\textwidth]{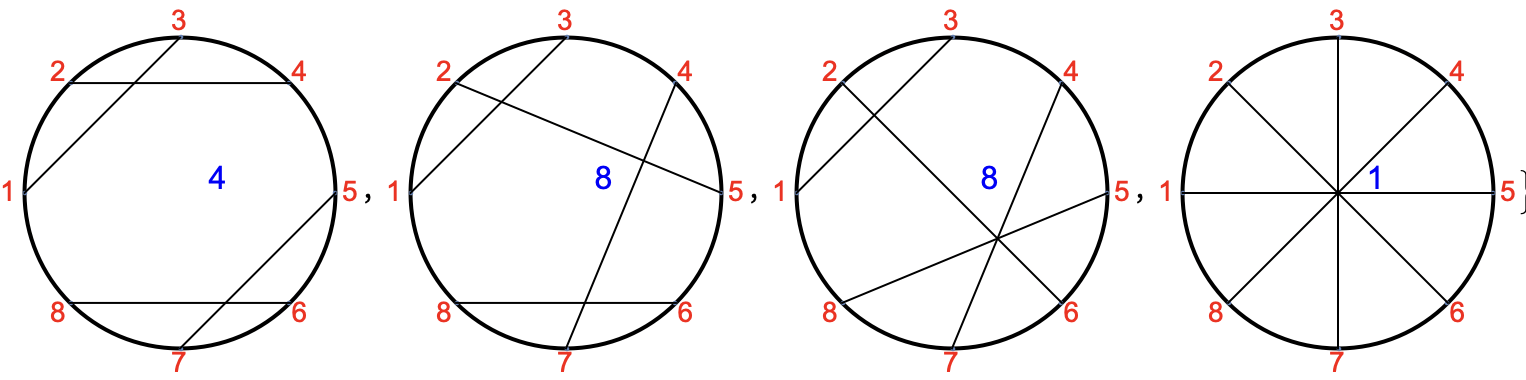}}
\end{center} 
\caption{\small{ The primitive diagrams of  order 8, type $[2^4]$ and genus 2, 
with their weight in blue}}
\label{Prim-HZ8}
\end{figure}

 \begin{figure}
\begin{center}
{\includegraphics[width=.42\textwidth]{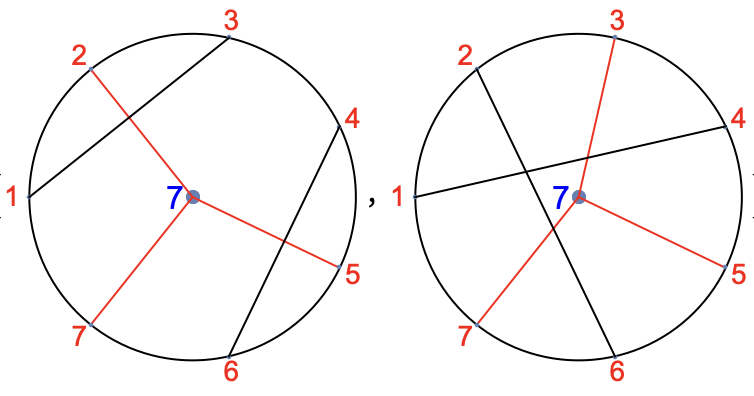}}
\end{center} 
\caption{\small{ The primitive diagrams of order 7, type $[2^2\,3]$ and genus 2, with the sum of weights equal to 14}}
\label{Prim7-322}
\end{figure}

 \begin{figure}
\begin{center}
{\includegraphics[width=.95\textwidth]{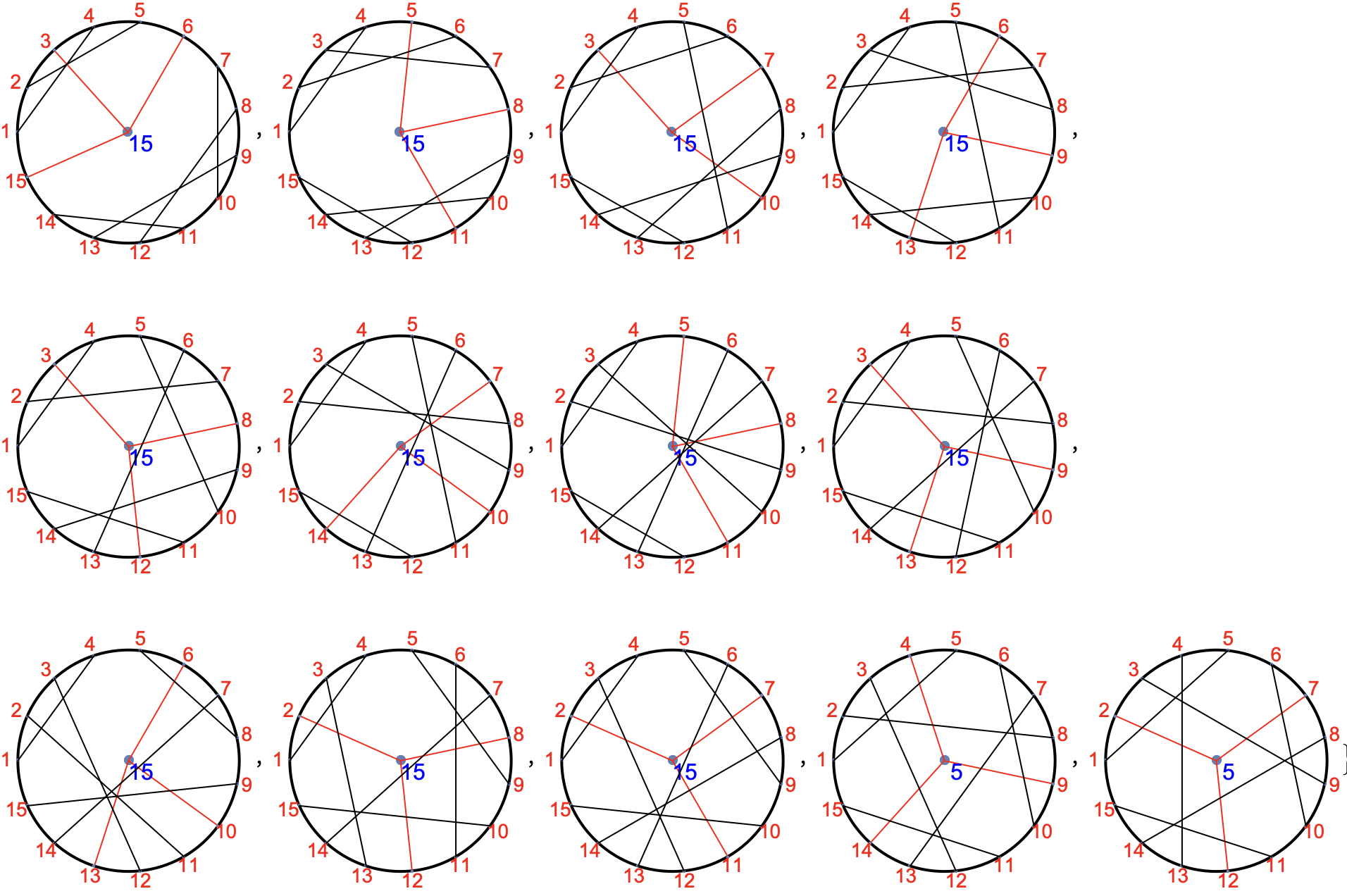}}\\
\end{center} 
\caption{\small{ The primitive diagrams of order 15, type $[2^6\,3]$  and genus 2, with the sum of weights equal to $175$}}
\label{Prim15-32222222}
\end{figure}

\begin{figure}
\begin{center}
{\includegraphics[width=.22\textwidth]{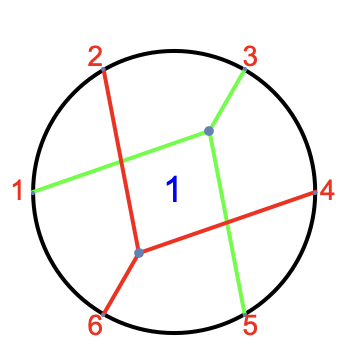}}
\end{center} 
\caption{\small{The primitive diagram of order 6, type $[3^2]$  and genus 2, of weight  1}}
\label{Prim6-33}
\end{figure}

\begin{figure}
\begin{center}
{\includegraphics[width=.6\textwidth]{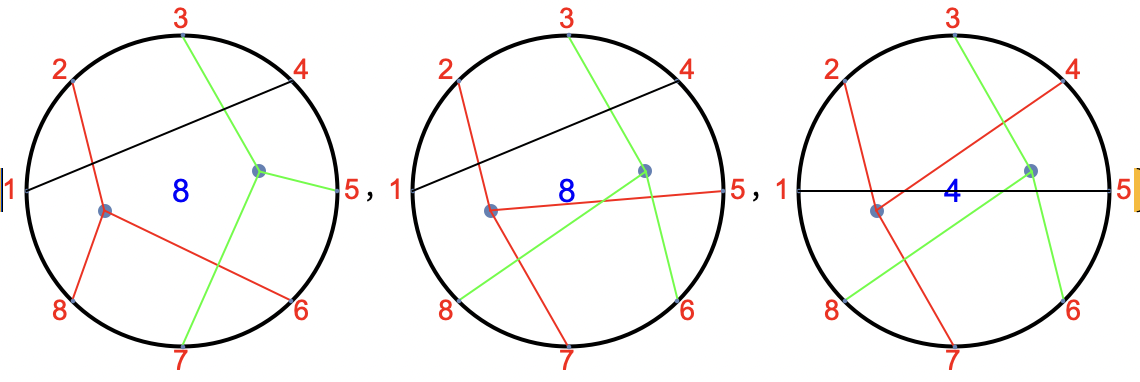}}
\end{center} 
\caption{\small{ The primitive diagrams of order 8, type $[2\,3^2]$  and genus 2,  of total weight  20}}
\label{Prim8-332}
\end{figure}

\begin{figure}
\begin{center}
{\includegraphics[width=.6\textwidth]{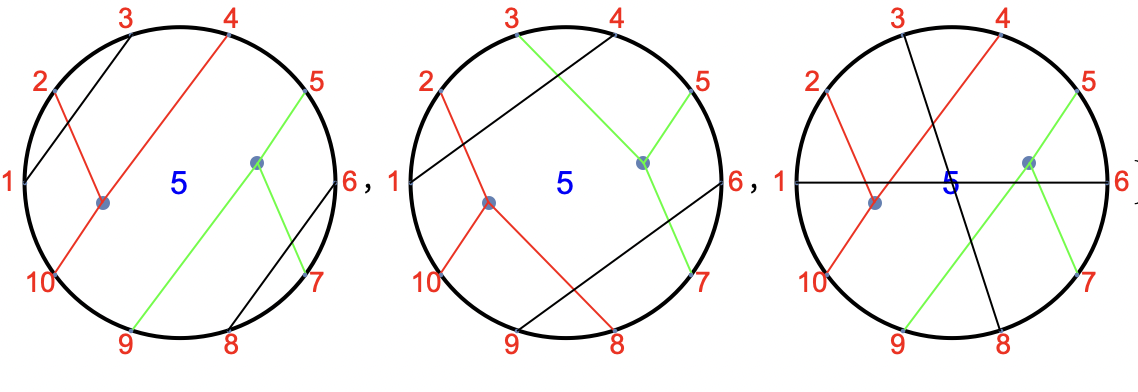}}
\end{center} 
\caption{\small{ The 3 semi-primitive diagrams of order 10, type $[2^2\, 3^2]$,  and genus 2, with the sum of weights equal to 15}}
\label{SemiPrim10-3322}
\end{figure} 

\begin{figure}
\begin{center}
{\includegraphics[width=.4\textwidth]{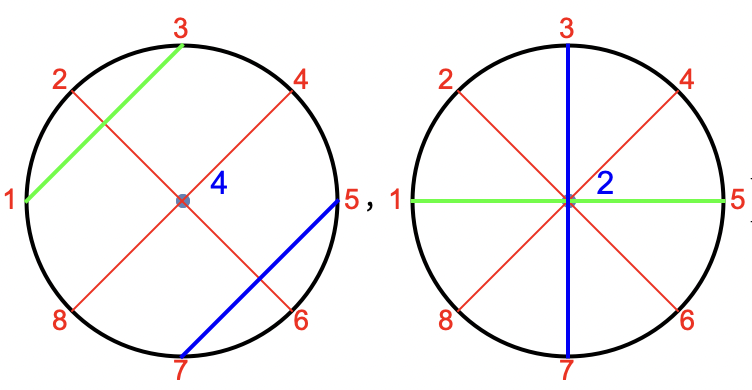}}
\end{center} 
\caption{\small{ The 2 primitive diagrams of order 8, type $[2^2\,4]$,  and genus 2, with the sum of weights equal to 6 }}
\label{Prim8-422}
\end{figure} 

\begin{figure}
\begin{center}
{\includegraphics[width=.4\textwidth]{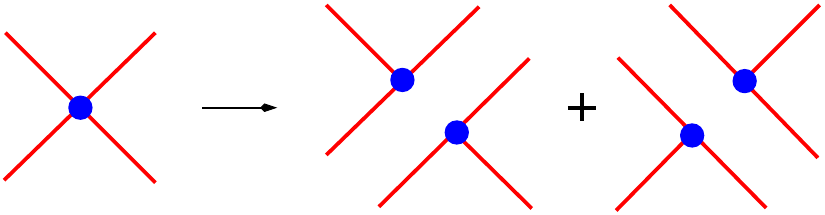}}
\end{center} 
\caption{\small{ The splitting procedure, by which here a 4-vertex is split into two 3-vertices}}
\label{splitting}
\end{figure} 


\subsection{Particular cases}
\subsubsection{Genus 2 partitions of $n=2p$ into $p$ doublets}
  In the simplest case where only $\kappa_2\ne 0$  (and set equal to 1 with no loss of generality),
the  primitive diagrams are of order $n\le 18$ -- a sample of which 
is shown in Fig. \ref{Prim-HZ8} \footnote{All genus 2 primitive and semi-primitive diagrams may be found on  \\ 
\url{https://www.lpthe.jussieu.fr/~zuber/Z_UnpubPart.html}}. 
They involve only 2-lines and their dressing is given by the expression (\ref{z2}) above. Thus
 \bea\nonumber  
 Z^{(2)}(x; \kappa_2=1, \kappa_{i\ne 2}= 0) 
 &=& \\    \nonumber  \frac{ \tY_2(x)}{(1-2 \kk2 x^2 Z^{(0)}(x))}\!\!\!\!\!\!\!\!   &&\!\!\!\! \!\!\!\! \Big(21 \tX_2^3 (x)+ 168 \tX_2^4(x)   +483 \tX_2^5(x) + 651   \tX_2^6(x) + 420  \tX_2^7(x) +   105  \tX_2^8(x) \Big)\eea
  with the notations of (\ref{deftX2tY2}). 
 After some substantial algebra (carried out by Mathematica), one finds
\be\label{HZg2} 
Z^{(2)}(x; \kappa_2=1, \kappa_{i\ne 2}= 0) = \frac{21  x^8 (1+\kk2 x^2)}{(1-4 \kk2 x^2)^{11/2}}\ee
 in agreement with the results of \cite{HZ}.  

 \subsubsection{Genus 2 partitions of $n=3p$ into $p$ triplets}
 We now assume as in sect. \ref{nto3p} that only $\kappa_3\ne 0$ (and equals 1 with no loss of generality).
  Let $s:=\sin \left(\frac{1}{3} \sin ^{-1}\left(\frac{3}{2} \sqrt{3} x^{3/2}\right)\right)$. Then, \newa{following (\ref{Z2genus2})}, 
   $Z^{(2)}$ takes the fairly cumbersome form
\bea\nonumber
Z^{(2)}(x; \kappa_3=1; \kappa_{i\ne 3}=0)&=& \\
\nonumber &&\hskip-5cm\frac{192 s^6 x^6 \left(8 s^3 \left(128 \left(11264 s^9+8676 \sqrt{3} s^6 x^{3/2}+3105 s^3 x^3\right)+9315 \sqrt{3}
   x^{9/2}\right)+729 x^6\right)}{
\left(2 \cos \left(\frac{1}{3} \Arccos\big(1-\frac{27 x^3}{2}\big)\right)-1\right)
\left(9 \sqrt{x^3}-4 \sqrt{3} \sin   \left(\frac{1}{3} \Arcsin\big(\frac{3 \sqrt{3x^3}}{2}\big)\right)\right)^{10} }\eea
(compare with the denominator of $Z^{(1)}$ in (\ref{FCg1}).  The first terms of the series expansion read
$$ x^6 + 144 x^9 + 6046 x^{12} + 149674 x^{15} + 2771028 x^{18} + 42679084 x^{21}+\cdots $$
One finds again a singular behaviour of the form
$$Z^{(2)}(x; \kappa_3=1; \kappa_{i\ne 3}=0) \sim \frac{\mathrm{const.}}{(\frac{4}{27}-x^3)^{11/2}}\,.$$

\subsubsection{Total number of genus 2 partitions}
Taking all $\kappa$'s equal to 1 (and possibly $\kappa_1=0$), as in sect. \ref{totalnbr01}, hence $W(x)=x/(1-x)$ or $\widehat W(x)=x^2/(1-x)$, 
we compute by (\ref{genus2}) the GF of the total number of genus 2 partitions (with or without singletons), 
and we recover the result of Cori and Hetyei \cite{CoriH17}
\be\label{Z2x} Z^{(2)}(x;\k_i=1) = \frac{x^6 (1 + 6 x - 19 x^2 + 21 x^3)}{(1 - 4 x)^{11/2}}\,, \ee
and  also  
\be\label{Z2xsf} Z^{(2)}(x;\k_1=0; \k_{i>1}=1)=\frac{x^6 (1 + 10 x + 5 x^2 + 5 x^3 + 9 x^4)}{(1 - 2 x - 3 x^2)^{11/2}}\ee
 in accordance with 
(\ref{ZZhat}).
\moi{\tiny Checked in FromGenus0ToGenus1ToGenus2.nb}

\subsubsection{Genus 2 partitions into $r$ parts}
The two-variable GF of the number of genus 2 partitions into a given number of parts is obtained as in sect. \ref{FGxy01}
by setting all $\kappa_i=y$. Theorem 2 leads to 
\bea \label{FGxy2}  Z^{(2)}(x,y) &=& \frac{x^6 y^2\, p(x,y)}{( (1 +x -x y)^2 -4x)^{11/2}} \\
\nonumber p(x,y) &=&  1- x(4-10 y) + x^2(6-10y-15y^2) -x^3(4+10y-39y^2+4y^3) \\
  \nonumber  &&\qquad+x^4(1+10y-15y^2-4y^3+8y^4)
\eea
as first derived by Cori--Hetyei~\cite{CoriH17}.  
Similar formulae are obtained if singletons are excluded
\bea \label{FGxy2sf}  \widehat Z^{(2)}(x,y) &=& \frac{x^6 y^2\, \hat p(x,y)}{( (1 -x )^2 -4x^2 y)^{11/2}} \\
\nonumber \hat p(x,y) &=&  1+ x(-4+14 y) + x^2(6-22y+21y^2) +x^3(-4+2y+7y^2) \\
  \nonumber  &&\qquad+x^4(1+6y-19y^2+21y^3)\,.
\eea  

The counting of genus 2 partitions into $r$ parts is then obtained by identifying the coefficient of $y^r$ in (\ref{FGxy2}). 
For example,  for $r=2$ (partitions into two parts with or without singleton)
\bea \nonumber Z^{(2)}(x; r=2)=\widehat Z^{(2)}(x; r=2)
&=&   \frac{x^6}{(1-x)^7} 
\\ \nonumber &=& \sum_{n\ge 6} {n\choose 6} x^n
 \\ \nonumber
 &=&  \sum_{n\ge 6} x^n \sum_{p=2}^{n-2} \frac{n}{6} {p-1\choose 2}{n-p-1\choose 2}   
\eea
in agreement with a general result for $r=2$ and arbitrary genus~\cite{comp}. 
For $r=3$ (partitions into three parts without singleton)
\be \nonumber Z^{(2)}(x; r=3)= \frac{14 x^7 (1 + 2 x)}{(1 -x)^9}=14 \sum_{n\ge 7} {n\choose 7} \frac{3n-13}{8} x^n \,.\ee


\section{Conclusion and perspectives}
 In principle the method could be extended to higher genus, but at the price of an increasing number
of (semi-)primitive diagrams, whose set remains to be listed, \newa{with an Ansatz of the form
\be\label{Ansatzg} Z^{(g)}(x) = \frac{\sum \mathrm{dressing\ of\ (semi-)primitive\ diagrams\ of\ genus}\ g }{1-\sum_n  n \kappa_n x^n (Z^{(0)}(x))^{ n-1} }\,. \ee}
For instance, in genus 3, primitive diagrams may occur up to  $n=30$
and they start at order $n=12$.
An Ansatz for partitions into doublets (\ie of type $[2^p]$), for $g=3$ is thus
$$ Z^{(3)}(x; \kappa_2=1, \kappa_{i\ne 2}= 0) = \frac{\tY_2(x) \tX_2^5(x)}{ (1-2 \kk2 x^2 Z^{(0)}(x))}  \sum_{j=0}^9  a_j \tX_2^j(x) $$
in which the numerical coefficients $a_j$ \newa{count the primitives of type $[2^{j+6}]$ and} may be determined against 
the known result of \cite{WL1, HZ}
  \be\label{HZg3} Z^{(3)}(x; \kappa_2=1, \kappa_{i\ne 2}= 0) =\frac{11x^{12} (135 + 558 x^2 + 158 x^4)}{(1 - 4 x^2)^{17/2}}\,. \ee
hence
\bea\nonumber Z^{(3)}(x; \kappa_2=1, \kappa_{i\ne 2}= 0) &=& 
\frac{11 \tY_2(x)    \tX_2^5(x)}{ (1-2 \kk2 x^2 Z^{(0)}(x))}  \, \Big(135   +  2313 \tX_2(x)   +  15728 \tX_2^2(x)   +57770    \tX_2^3(x) 
   \\ \nonumber  + 128985 \tX_2^4(x) \!\!\!\! &+&\!\!\!\! 183955\tX_2 ^5(x)  + 169078 \tX_2^6(x)   + 97188  \tX_2^7(x)  +  31850 \tX_2^8(x)   +4550  \tX_2^9 (x)  \Big)\,.
  \eea
\omit{\red [while for higher genus, we may express the known exact results \cite{WL1,WL2,HZ} as
\bea \nonumber Z^{(2)}(x; \kappa_2=1, \kappa_{i\ne 2}= 0) &=&\frac{21 x^8 (1 + x^2)}{(1 - 4 x^2)^{11/2}}\\
\label{highergenusHZ} Z^{(3)}(x; \kappa_2=1, \kappa_{i\ne 2}= 0) &=&\frac{
11x^{12} (135 + 558 x^2 + 158 x^4)}{(1 - 4 x^2)^{17/2}}\\
 \nonumber Z^{(4)}(x; \kappa_2=1, \kappa_{i\ne 2}= 0) &=&\frac{
 143 x^{16} (1575 + 13689 x^2 + 18378 x^4 + 2339 x^6)}{(1 - 4 x^2)^{23/2}}
  \eea etc. ]\normalcolor}

{Likewise, in genus 4,
\bea\label{HZg4}  Z^{(4)}(x; \kappa_2=1, \kappa_{i\ne 2}= 0) &=&    \frac{143 x^{16}(1575 + 13689 x^2 + 18378 x^4 + 
    2339 x^6)}{(1 - 4 x^2)^{23/2}}\\  \nonumber
&=&\frac{143 \tY_2 \tX_2^7}{(1-2 \kk2 x^2 Z^{(0)}(x))}
\Big(1575 + 43614 \tX_2 + 497277 \tX_2^2 + 3194702 \tX_2^3 + 
   13162499 \tX_2^4 \\   
  \nonumber   \!\!\!\! &+&\!\!\!\!  37212840 \tX_2^5   +  74956749 \tX_2^6 + 
   109645557 \tX_2^7 + 117063972 \tX_2^8 + 90449979 \tX_2^9  \\
  \nonumber   \!\!\!\! &+&\!\!\!\! 
   49312410 \tX_2^{10} + 18008865 \tX_2^{11} + 3956750 \tX_2^{12} + 
   395675 \tX_2^{13}\Big)\,.
  \eea}

We end this paper with \new{a few remarks on} some intriguing issues.  \\
There is some evidence of a  universal singular behaviour of all generating functions, 
\be\label{singbehr}Z^{(g)}(x)\sim (x_0-x)^{\oh -3g}\ee 
as can be  seen on the partitions into doublets (\ref{HZg1},\ref{HZg2},\ref{HZg3}, \ref{HZg4} ),  and for $ g=1,2$ on other cases.
This would imply a large $n$ behaviour of coefficients $C_{n, [\alpha]}^{(g)}$ (for  appropriately rescaled patterns $\alpha $)
of the form
$$ C_{n, [\alpha]}^{(g)}  \sim \text{const} \  x_0^{-n-3g+\oh} \,  n^{3g-\oh}    \qquad \mathrm{as} \ n, [\alpha] \mathrm{\ grow\ large} \,.$$ 
This type of singularity of the GF  and the associated asymptotic behaviour have been observed 
in the parallel problem of enumeration of unicellular maps by Chapuy~\cite{Chapuy}, who interpreted  
 the number $6g-1$ as the number of edges in his dominant ``schemes"  (the analogues of our primitives).
That the same behaviour appears in the present context of partitions indicates that the restriction of maps
due to the restricted crossing constraint discussed in sect. \ref{matint} is ``irrelevant" (in the sense of critical phenomena), \ie does not affect the singular behaviour. 
 The ``critical exponent" $\oh -3g$
is also familiar to physicists in the context of boundary loop models and Wilson loops~\cite{Kos}. \new{Such a connection is natural in the case 
of partitions into doublets, since it is known that in that case, the counting amounts to computing the expectation value of $\tr M^n$
in a Gaussian matrix integral, hence for large $n$, of a large loop. That the same singular or asymptotic behaviour takes place 
in (all ?) other cases seems to indicate that an effective Gaussian theory takes place in that limit.}\footnote{I'm grateful to Ivan Kostov for discussions on that point}
\\[6pt]
A natural question is whether  the Topological Recurrence of Chekhov, Eynard  and Orantin \cite{EyOCh} is relevant for the counting of partitions and is 
related to or independent of the approach of this paper.
\\[6pt]
\omit{The  dressing procedure of primitive diagrams  is reminiscent of techniques used in QFT, under the name of Dyson-Schwinger equation. 
This leads to another natural question: is there an
effective theory --matrix-like or otherwise-- that would lie behind the screen ?
\\[6pt]}
As mentioned in the introduction, the  formulae derived in this paper yield an interpolation between 
expansions on ordinary and  on free cumulants. What is the relevance of this interpolation? How does it compare with other existing 
interpolations ?
\\[10pt]}
\normalcolor 
All these questions are left for future investigation.


\bigskip
\noindent {\bf Acknowledgements.} \neww{It is a pleasure to thank Guillaume Chapuy, 
 Philippe Di Francesco,  Elba Garcia-Failde and Ivan Kostov for discussions and comments, \newa{and Colin McSwiggen for suggesting amendments of this paper}. I'm particularly grateful to Robert Coquereaux for a  careful reading of a first draft of the manuscript  and for providing me with very efficient Mathematica codes. }

\normalcolor

\end{document}